\documentclass[12pt]{amsart}
\usepackage{amssymb}
\usepackage{amsfonts}
\usepackage{amscd}
\usepackage[all]{xypic}

\swapnumbers

\def\limproj{\mathop{\oalign{lim\cr
\hidewidth$\longleftarrow$\hidewidth\cr}}}

\def\ac{{\rm ac}}
\def\pr{{\rm pr}}

\def\Var{{\rm Var}}
\def\hsp{{\rm hsp}}

\def\Sp{{\rm Sp}}

\def\HS{{\rm HS}}

\def\GM{{\mathbf G}_m}

\def\AA{{\mathbf A}}

\def\CC{{\mathbf C}}

\def\GG{{\mathbf G}}

\def\LL{{\mathbf L}}

\def\NN{{\mathbf N}}

\def\PP{{\mathbf P}}
\def\QQ{{\mathbf Q}}
\def\RR{{\mathbf R}}

\def\ZZ{{\mathbf Z}}

\def\cA{{\mathcal A}}

\def\cI{{\mathcal I}}
\def\cJ{{\mathcal J}}

\def\cL{{\mathcal L}}
\def\cM{{\mathcal M}}

\def\cO{{\mathcal O}}

\def\cS{{\mathcal S}}

\def\cW{{\mathcal W}}
\def\cX{{\mathcal X}}
\def\cY{{\mathcal Y}}

\mathchardef\alphag="7C0B \mathchardef\betag="7C0C
\mathchardef\gammag="7C0D \mathchardef\deltag="7C0E
\mathchardef\varepsilong="7C22 \mathchardef\varphig="7C27
\mathchardef\psig="7C20 \mathchardef\zetag="7C10
\mathchardef\epsilong="7C0F \mathchardef\rhog="7C1A
\mathchardef\taug="7C1C \mathchardef\upsilong="7C1D
\mathchardef\iotag="7C13 \mathchardef\thetag="7C12
\mathchardef\pig="7C19 \mathchardef\sigmag="7C1B
\mathchardef\etag="7C11 \mathchardef\omegag="7C21
\mathchardef\kappag="7C14 \mathchardef\lambdag="7C15
\mathchardef\mug="7C16 \mathchardef\xig="7C18
\mathchardef\chig="7C1F \mathchardef\nug="7C17
\mathchardef\varthetag="7C23 \mathchardef\varpig="7C24
\mathchardef\varrhog="7C25 \mathchardef\varsigmag="7C26
\mathchardef\Omegag="7C0A \mathchardef\Thetag="7C02
\mathchardef\Sigmag="7C06 \mathchardef\Deltag="7C01
\mathchardef\Phig="7C08 \mathchardef\Gammag="7C00
\mathchardef\Psig="7C09 \mathchardef\Lambdag="7C03
\mathchardef\Xig="7C04 \mathchardef\Pig="7C05
\mathchardef\Upsilong="7C07

\newtheorem{theorem}[subsection]{Theorem}
\newtheorem{lem}[subsection]{Lemma}
\newtheorem{cor}[subsection]{Corollary}
\newtheorem{prop}[subsection]{Proposition}

\newtheorem{conj}[subsection]{Conjecture}

\theoremstyle{definition}

\newtheorem{def-prop}[subsubsection]{Definition-Proposition}

\theoremstyle{remark}
\newtheorem{remark}[subsection]{Remark}

\theoremstyle{plain}

\numberwithin{equation}{subsection}

\def\boxit#1#2{\setbox1=\hbox{\kern#1{#2}\kern#1}%
\dimen1=\ht1 \advance\dimen1 by #1 \dimen2=\dp1
\advance\dimen2 by #1
\setbox1=\hbox{\vrule height\dimen1 depth\dimen2\box1\vrule}%
\setbox1=\vbox{\hrule\box1\hrule}%
\advance\dimen1 by .4pt \ht1=\dimen1 \advance\dimen2 by
.4pt \dp1=\dimen2 \box1\relax}

\newcommand{\sur}[2]{\genfrac{}{}{0pt}{}{#1}{#2}}

\def\AA{{\mathbf A}}

\def\CC{{\mathbf C}}

\def\GG{{\mathbf G}}

\def\LL{{\mathbf L}}

\def\NN{{\mathbf N}}

\def\PP{{\mathbf P}}
\def\QQ{{\mathbf Q}}
\def\RR{{\mathbf R}}

\def\ZZ{{\mathbf Z}}

\def\cA{{\mathcal A}}

\def\cI{{\mathcal I}}
\def\cJ{{\mathcal J}}

\def\cL{{\mathcal L}}
\def\cM{{\mathcal M}}

\def\cO{{\mathcal O}}

\def\cS{{\mathcal S}}

\def\cW{{\mathcal W}}
\def\cX{{\mathcal X}}
\def\cY{{\mathcal Y}}

\mathchardef\alphag="7C0B \mathchardef\betag="7C0C
\mathchardef\gammag="7C0D \mathchardef\deltag="7C0E
\mathchardef\varepsilong="7C22 \mathchardef\varphig="7C27
\mathchardef\psig="7C20 \mathchardef\zetag="7C10
\mathchardef\epsilong="7C0F \mathchardef\rhog="7C1A
\mathchardef\taug="7C1C \mathchardef\upsilong="7C1D
\mathchardef\iotag="7C13 \mathchardef\thetag="7C12
\mathchardef\pig="7C19 \mathchardef\sigmag="7C1B
\mathchardef\etag="7C11 \mathchardef\omegag="7C21
\mathchardef\kappag="7C14 \mathchardef\lambdag="7C15
\mathchardef\mug="7C16 \mathchardef\xig="7C18
\mathchardef\chig="7C1F \mathchardef\nug="7C17
\mathchardef\varthetag="7C23 \mathchardef\varpig="7C24
\mathchardef\varrhog="7C25 \mathchardef\varsigmag="7C26
\mathchardef\Omegag="7C0A \mathchardef\Thetag="7C02
\mathchardef\Sigmag="7C06 \mathchardef\Deltag="7C01
\mathchardef\Phig="7C08 \mathchardef\Gammag="7C00
\mathchardef\Psig="7C09 \mathchardef\Lambdag="7C03
\mathchardef\Xig="7C04 \mathchardef\Pig="7C05
\mathchardef\Upsilong="7C07

\DeclareMathOperator*{\Spec}{Spec}
\def\ord{{\rm ord}}

\begin{document}

\title[Iterated vanishing cycles]{Iterated vanishing cycles,
convolution,  and a motivic analogue of a conjecture of
Steenbrink}

\author{Gil Guibert}

\address{39 quai du Halage,
94000 Cr\'eteil, France}

\email{guibert9@wanadoo.fr}

\author{Fran\c cois Loeser}

\address{{\'E}cole Normale Sup{\'e}rieure,
D{\'e}partement de math{\'e}matiques et applications, 45
rue d'Ulm, 75230 Paris Cedex 05, France (UMR 8553 du CNRS)}
\email{Francois.Loeser@ens.fr}
\urladdr{http://www.dma.ens.fr/~loeser/}

\author{Michel Merle}

\address{Laboratoire J.-A. Dieudonn\'e,
Universit\'e de Nice - Sophia Antipolis, Parc Valrose,
06108 Nice Cedex 02, France (UMR 6621 du CNRS)}
\email{merle@math.unice.fr}
\urladdr{http://www-math.unice.fr/membres/merle.html}

\maketitle
\section{Introduction}Let us start by recalling the statement of
Steenbrink's conjecture.
Let $f : X \rightarrow \AA^1$ be a function on a smooth
complex algebraic variety.
Let $x$ be a closed point of $f^{-1} (0)$.
Steenbrink introduced \cite{St1} the notion of
the spectrum $\Sp (f, x)$ of $f$ at $x$.
It is a fractional Laurent polynomial
$\sum_{\alpha \in \QQ} n_{\alpha} t^{\alpha}$, $n_{\alpha}$ in
$\ZZ$, which
is constructed using the action
of the monodromy on the mixed Hodge structure
on the cohomology of the Milnor fiber at $x$.
When $f$ has an isolated singularity at $x$,
all $n_{\alpha}$ are in $\NN$, and the exponents of $f$,
counted with multiplicity $n_{\alpha}$, are
exactly the rational numbers
$\alpha$ with $n_{\alpha}$ not zero.

Let us assume now that the singular locus of  $f$ is a curve
$\Gamma$, having $r$ local components $\Gamma_{\ell}$,
$1\leq \ell \leq r$, in a neighborhood of
$x$.
We denote by $m_{\ell}$ the multiplicity of $\Gamma_{\ell}$.
Let $g$ be a generic linear form vanishing at $x$
 (that is,  a function $g$
vanishing at $x$ whose differential at $x$ is a generic linear form).
For $N$ large enough, the function  $f+g^N$
has an isolated singularity at $x$.
In a neighborhood of
the complement $\Gamma_{\ell}^\circ$
to $\{x\}$
in $\Gamma_{\ell}$, we may view
$f$ as a family
of isolated hypersurface singularities
parametrized by $\Gamma_{\ell}^\circ$.
The cohomology of the Milnor fiber
of this hypersurface singularity
is naturally endowed with the action of
two commuting monodromies: the monodromy of the function and
the monodromy of a generator of the local fundamental group
of $\Gamma_{\ell}^\circ$.
We denote by $\alpha_{\ell, j}$ the exponents of
that isolated hypersurface singularity
and by
$\beta_{\ell, j}$ the corresponding rational
numbers in $[0, 1)$ such
that the complex numbers $\exp (2 \pi i \beta_{\ell, j})$
are the eigenvalues of the monodromy along
$\Gamma_{\ell}^\circ$.

\begin{conj}[Steenbrink \cite{Stconj}]\label{steenconj}For $N \gg 0$,
\begin{equation}\label{ooo}
\Sp (f + g^N, x)
-
\Sp (f, x) = \sum_{\ell, j}
t^{\alpha_{\ell, j} + (\beta_{\ell, j}/m_{\ell} N)}
\frac{1- t}{1 -t^{1/m_{\ell} N}}.
\end{equation}
\end{conj}

The conjecture of Steenbrink has been proved by M. Saito in
\cite{Sa3}, using his theory of mixed Hodge modules
\cite{Sa1}\cite{Sa2}. Later, A. N\'emethi and J. Steenbrink
\cite{ns} gave another proof, still relying on the theory
of mixed Hodge modules. Also, forgetting the integer part
of the exponents of the spectrum, (\ref{ooo}) has been
proved by D. Siersma \cite{Siersma} in terms of zeta
functions of the monodromy. Notice that, taking ordinary
Euler characteristics, (\ref{ooo}) specializes to a result
of I. Iomdin \cite{iomdin} who was the first to compare
vanishing cohomologies of $f$ and $f+g^N$.
The convention we use here, cf. (\ref{defhod}), to define
$\Sp (f, x)$ slightly differs from the original one and corresponds
to what is denoted by $\Sp' (f, x)$
in \cite{Sa3}.

Recently, using motivic integration, Denef and Loeser introduced
the motivic Milnor fiber $\cS_{f, x}$.
It is a virtual variety endowed
with 
an action of the group scheme $\hat \mu$ of roots of unity
 and the Hodge spectrum
$\Sp (f, x)$ can be retrieved from $\cS_{f, x}$, cf. \cite{barc}.
They also showed
that an analogue of the Thom-Sebastiani Theorem holds
for the motivic Milnor fiber.
This result was first stated in a (completed) Grothendieck ring \cite{ts}
of Chow motives and then extended to a
Grothendieck ring of virtual varieties endowed
with a $\hat \mu$-action
in \cite{looi} and \cite{barc}, using
a convolution product $\ast$ introduced in \cite{looi}.
It is also convenient to
slightly modify
the virtual varieties $\cS_{f, x}$,
which correspond to nearby cycles,
into virtual varieties
$\cS_{f, x}^\phi$
corresponding to vanishing cycles.

It is then quite natural to ask
for a motivic analogue of Steenbrink's conjecture
in terms
of motivic Milnor fibers. The present paper is devoted to give
a complete answer to that question.
Our main result, Theorem \ref{MT},
expresses (in its local version \ref{mc}), for $x$
a closed point where $f$ and $g$ both vanish and for $N \gg 0$, the difference
$\cS_{f, x}^\phi - \cS_{f+g^N, x}^\phi$ as
$\Psi_{\Sigma} (\cS_{g^N, x}(\cS_f^\phi))$,  where
$\cS_{g^N, x}(\cS_f^\phi)$ corresponds to iterated
motivic vanishing cycles
and $\Psi_{\Sigma}$ is a generalization
of the convolution product $\ast$.
In fact,  in Theorem \ref{MT}, we no longer
assume any
condition on the singular locus of $f$; also
$g$ is not assumed anymore to be a generic linear form and can be any function vanishing at $x$.
Formula (\ref{ooo}) may be deduced from
Theorem \ref{MT} by considering the Hodge spectrum.

\medskip

The plan of the paper is the following.
In section \ref{gr} we introduce the basic Grothendieck rings
we shall use.
Then, in section \ref{mvc}, we
recall the definition of the Motivic Milnor fiber and we
extend it to the whole Grothendieck ring.
Such an extension has also been done by F. Bittner
in \cite{bittner2}, using the weak factorization
Theorem and her work \cite{bittner} ; the construction we present here,
based on  motivic integration,
is quite different. We then extend the construction to the
equivariant setting, in order
to define
iterated vanishing cycles in the motivic framework in section \ref{ivcsec}.
In section \ref{cmr} we first define
our generalized convolution operator
$\Psi_{\Sigma}$ and explain its relation with the convolution
product $\ast$.
This gives us the opportunity to prove the associativity of
the convolution product $\ast$, a fact already mentioned
in \cite{barc}.
Then comes the heart of the paper, that is the proof of
Theorem \ref{MT}. We conclude the section
by explaining how one recovers the motivic Thom-Sebastiani
Theorem of \cite{ts}, \cite{looi} and
\cite{barc} from Theorem \ref{MT}.
The final section \ref{ssc} is devoted
to applications to the Hodge-Steenbrink spectrum,
in particular, we deduce Steenbrink's conjecture \ref{steenconj}
from
Theorem \ref{MT}.

\subsection*{}

We heartfully thank the referee for   
reading several versions of the paper with meticulous care
and for providing us
comments that were most helpful in improving the exposition and in eliminating many inaccuracies.
We also thank  A. Fernandez de Bobadilla and A. Melle for their comments.

\section{Grothendieck rings}\label{gr}
\subsection{}By a variety over of field $k$,
we mean a separated and reduced scheme of finite type over $k$.
If $X$ is a scheme, we denote by $\vert X\vert$ the corresponding reduced
scheme.
If an algebraic group $G$ acts on a variety $X$, we say
the action is good if every $G$-orbit
is contained in an affine open subset of $X$.
Let $Y$ be a variety over $k$ and let $p : A \rightarrow Y$ be an affine bundle for the Zariski topology
(the fibers of $p$ are affine spaces and the transition morphisms between
trivializing charts are affine). In particular the fibers of $p$ have the structure of affine spaces.
Let $G$ be a linear  algebraic group. A good action of $G$ on $A$ is said to be affine if it is
a lifting of a 
good action
on $Y$ and its restriction to all fibers is affine.
Note that affine actions on an 
affine bundle  extend to its relative projective bundle  compactification.

If $G$ is finite and  $X$ and $Y$ are two varieties
with good $G$-action,
we denote by $X \times^G Y$ the quotient
of the product
$X \times Y$ by the equivalence relation
$(gx, y) \equiv
(x, gy)$.
The action of $G$ on, say, the first factor of
$X \times Y$ induces a good $G$-action on $X \times^G Y$.

For $n \geq 1$, we denote by $\mu_n$ the group scheme
of $n$-th roots of unity
and
by $\hat \mu$ the projective limit
$\limproj \mu_n$ of
the projective system with transition morphisms
$\mu_{nd} \rightarrow \mu_n$ given by 
$x \mapsto x^d$. In this paper all $\hat \mu$-actions, and more
generally all $\hat \mu^r$-actions, will be assumed
to factorize through a finite quotient.

\subsection{}Throughout the paper
$k$ will be a field of characteristic zero.
For $S$ a variety over $k$, we denote by
$K_0 (\Var_S)$ the Grothendieck ring
of varieties over $S$, cf. \cite{barc}.
Let us recall it is generated by
classes of morphisms of varieties
$X \rightarrow S$, and that it is also generated by classes
of such morphisms with $X$ smooth over $k$ and it suffices to consider
relations for smooth varieties.
We denote by $\LL = \LL_S$ the class of the trivial
line bundle over $S$ and
set
$\cM_S$ for the localization
$K_0 (\Var_S) [\LL^{-1}]$.
As in \cite{lef}, let us consider Grothendieck rings
of varieties with $\hat \mu$-action.
They are defined similarly, using
the category
$\Var_S^{\hat \mu}$ of
varieties with good $\hat \mu$-action over $S$,
but adding the additional relation
\begin{equation}\label{addrelmu}
[Y \times \AA^n_k, \sigma] = [Y \times \AA^n_k, \sigma']
\end{equation}
if $\sigma$ and $\sigma'$ are two liftings of the
same $\hat \mu$-action on $Y$ to an 
affine 
action on 
$Y \times \AA^n_k$
We shall denote them by
$K_0 (\Var_S^{\hat \mu})$
and $\cM_S^{\hat \mu}$.
One can more generally replace
$\hat \mu$ by $\hat \mu^r$ in these definitions
and define
$K_0 (\Var_S^{\hat \mu^r})$
and $\cM_S^{\hat \mu^r}$.
In \cite{bittner2} Bittner considers similar
equivariant rings, but with an additional relation
a priori coarser than the one we use here.

\subsection{}\label{676}In the present paper, instead of varieties with $\hat \mu$-action
over $S$, we choose to  work
in the equivalent setting of varieties with
$\GM$-action with some additional structure.

Let $Y$ be a variety with good $\GM^r$-action. We say a
morphism $\pi : Y \rightarrow \GM^r$ is diagonally
monomial of weight
${\bf n}$ in $\NN_{>0}^r$, if $\pi (\lambda x) = \lambda^{\bf n} \pi (x)$ for all $\lambda$ in $\GM^r$ and
$x$ in $Y$. Fix ${\bf n}$ in $\NN_{>0}^r$. We denote by
$\Var_{S \times \GM^r}^{\GM^r, {\bf n}}$ the category of
varieties $Y \rightarrow S \times \GM^r$ over $S \times
\GM^r$ with good $\GM^r$-action such that furthermore, the
fibers of the projection $\pi_1 : Y \rightarrow S$ are
$\GM^r$-invariant and the projection $\pi_2 :  Y
\rightarrow \GM^r$ is diagonally monomial of weight ${\bf n}$. We
define the Grothendieck group $K_0 (\Var_{S \times
\GM^r}^{\GM^r, {\bf n}})$ as the free abelian group on
isomorphism classes of objects  $Y \rightarrow S \times
\GM^r$ in $\Var_{S \times \GM^r}^{\GM^r , {\bf n}}$, modulo
the relations
\begin{equation}[Y \rightarrow S \times \GM^r]  =  [Y'
\rightarrow S \times \GM^r] + [Y \setminus Y' \rightarrow S
\times \GM^r]
\end{equation} for $Y'$ closed $\GM^r$-invariant in $Y$
and, for $f : Y \rightarrow S \times \GM^r$ in $\Var_{S
\times \GM^r}^{\GM^r,  {\bf n}}$, 
\begin{equation}\label{addrel}[Y \times \AA^n_k
\rightarrow S \times \GM^r, \sigma] = [Y \times \AA^n_k
\rightarrow S \times \GM^r, \sigma']
\end{equation}if $\sigma$ and
$\sigma'$ are two liftings of the same $\GM^r$-action on
$Y$ to 
affine actions, 
the morphism $Y \times \AA^n_k  \rightarrow S \times
\GM^r$ being the composition of $f$ with the projection on the
first factor. 
Of course in (\ref{addrel}), instead of the trivial affine
bundle we could have considered
any affine bundle over $Y$.

Fiber product over $S\times \GM^r$ with diagonal action
induces a product in the category $\Var_{S \times
\GM^r}^{\GM^r, {\bf n}}$, which allows to endow $K_0
(\Var_{S \times \GM^r}^{\GM, {\bf n}})$ with a natural ring
structure. Note that the unit $1_{S\times \GM^r}$ for the
product is the class of the identity morphism  on $S \times
\GM^r$, the  $\GM^r$-action on $S \times \GM^r$ being the
trivial one on $S$ and standard multiplicative translation
on $\GM^r$.There is a natural structure of $K_0
(\Var_k)$-module on $K_0 (\Var_{S \times \GM^r}^{\GM^r,
{\bf n}})$. We denote by $\LL_{S\times \GM^r} = \LL$ the
element $\LL \cdot 1_{S\times \GM^r}$ in this module, and
we set $\cM_{S \times \GM^r}^{\GM^r, {\bf n}} = K_0
(\Var_{S \times \GM^r}^{\GM^r, {\bf n}}) [\LL^{-1}]$.

If $f : S \rightarrow S'$ is a morphism of varieties, composition with $f$ leads to a
push-forward morphism
$f_! : \cM_{S \times \GM^r}^{\GM^r, {\bf n}} \rightarrow \cM_{S' \times \GM^r}^{\GM^r, {\bf n}}$,
while fiber product
leads to a
pull-back morphism
$f^*: \cM_{S' \times \GM^r}^{\GM^r, {\bf n}} \rightarrow \cM_{S \times \GM^r}^{\GM^r, {\bf n}}$
(these morphisms may already be defined at the $K_0$-level).

\subsection{}For  ${\bf n}$ in
$\NN_{>0}^r$, we denote by $\mu_{\bf n}$
the group $\mu_{n_1} \times \cdots \times \mu_{n_r}$.
We consider the functor
\begin{equation}
G_{\bf n} :
\Var_{S \times \GM^r}^{\GM^r, {\bf n}}
\longrightarrow
\Var_S^{\mu_{\bf n}}
\end{equation}
assigning to $p : Y \rightarrow S \times \GM^r$
the fiber at $1$ of the morphism
$Y \rightarrow \GM^r$ obtained by composition with
projection on the second factor. Note that
this fiber carries a natural $\mu_{\bf n}$-action by the monomiality
assumption.

On the other side, if
$f: X \rightarrow  S$ is  a variety over $S$ with good $\mu_{\bf n}$-action,
we may
consider the variety
$F_{\bf n} (X) := X \times^{\mu_{\bf n}} \GM^r$
and view it as a variety over
$S \times \GM^r$
by sending the class of $(x, \lambda)$ to
$(f (x), \lambda^{\bf n})$.
The standard $\GM^r$-action by multiplicative translation
on $\GM^r$ induces a $\GM^r$-action on $F_{\bf n} (X)$.
Note that the second projection is  diagonally monomial
of weight
${\bf n}$, hence $F_{\bf n}$ is in fact a functor
\begin{equation}
F_{\bf n} : \Var_S^{\mu_{\bf n}} \longrightarrow
\Var_{S \times \GM^r}^{\GM^r, {\bf n}}.
\end{equation}

\begin{lem}
The functors
$F_{\bf n}$ and $G_{\bf n}$
are mutually quasi-inverse, so that
the  categories
$\Var_S^{\mu_{\bf n}}$ and
$\Var_{S \times \GM^r}^{\GM^r, {\bf n}}$
are equivalent.
\end{lem}

\begin{proof}It is quite clear that
$G_{\bf n} (F_{\bf n} (X))$ is isomorphic to $X$, for $X$ in
$\Var_S^{\mu_{\bf n}}$.
For $X$ in $\Var_{S \times \GM^r}^{\GM^r, {\bf n}}$, set $Y := G_{\bf n} (X)$.
We have a natural morphism
$Y \times \GM^r \rightarrow X$
sending $(y, \lambda)$ to $\lambda y$.
Clearly this morphism induces an isomorphism
between
$Y \times ^{\mu_{\bf n}} \GM^r$ and $X$
in
$\Var_{S \times \GM^r}^{\GM^r, {\bf n}}$.
\end{proof}

We consider the partial order ${\bf n} \prec {\bf m}$ on
$\NN_{>0}^r$
given by divisibility of each coordinates, that is,
${\bf n} \prec {\bf m}$ if ${\bf n} = {\bf k} {\bf m}$
for some $\bf k$ in $\NN_{>0}^r$.
If
${\bf n} = {\bf k} {\bf m}$, we have a natural functor
\begin{equation}\label{dlel}\theta^{\bf m}_{\bf n} : \Var_{S \times \GM^r}^{\GM^r, {\bf m}} \longrightarrow
\Var_{S \times \GM^r}^{\GM^r, {\bf n}},
\end{equation}
sending
$X \rightarrow S \times \GM^r$
to the same object, but with the action $\lambda \mapsto \lambda x$
on $X$ replaced by $\lambda \mapsto \lambda^{\bf k} x$.
We define the category
$\Var_{S \times \GM^r}^{\GM^r}$ as the colimit of the inductive system
of categories
$\Var_{S \times \GM^r}^{\GM^r, {\bf n}}$.
We define
$K_0 (\Var_{S \times \GM^r}^{\GM^r})$
and $\cM_{S \times \GM^r}^{\GM^r}$
as in \ref{676}.
Clearly, $K_0 (\Var_{S \times \GM^r}^{\GM^r})$
and $\cM_{S \times \GM^r}^{\GM^r}$
are respectively the colimits of the rings
$K_0 (\Var_{S \times \GM^r}^{\GM^r, {\bf n}})$
and $\cM_{S \times \GM^r}^{\GM^r, {\bf n}}$.
Since the category
$\Var_S^{{\hat \mu}^r}$ is the colimit of the categories
$\Var_S^{{\hat \mu_{\bf n}}}$, we have the following statement:

\begin{prop}\label{equiv}There is a unique pair of functors
\begin{equation}
G :
\Var_{S \times \GM^r}^{\GM^r}
\longrightarrow
\Var_S^{{\hat \mu}^r}
\end{equation}
and
\begin{equation}
F : \Var_S^{{\hat \mu}^r}
\longrightarrow
\Var_{S \times \GM^r}^{\GM^r}
\end{equation}
that restrict to $G_{\bf n}$ and $F_{\bf n}$
for every ${\bf n}$.
They are mutually quasi-inverse.
In particular $G$ induces canonical
isomorphisms
 \begin{equation}\label{comp}
K_0 (\Var_{S \times \GM^r}^{\GM^r})
\simeq
K_0 (\Var_S^{{\hat \mu}^r})
\quad
\text{and}
\quad
\cM_{S \times \GM^r}^{\GM^r}
\simeq
\cM_S^{{\hat \mu}^r}
\end{equation}
compatible with the operations $f_!$ and $f^*$. \qed
\end{prop}

\subsection{}\label{twist}Let $Y$ be a variety with good $\GM^r$-action. We say a
morphism $\pi : Y \rightarrow \GM^r$ is
monomial 
if it is equivariant with respect to some
transitive $\GM^r$-action on $\GM^r$
(see section \ref{resol} for monomial morphisms that are
not diagonally monomial morphisms).
More generally
consider a 
variety $(p, \pi): Y \rightarrow S \times \GM^r$ over $S \times
\GM^r$ with good $\GM^r$-action such that furthermore, the
fibers of $p : Y \rightarrow S$ are
$\GM^r$-invariant and $\pi:  Y
\rightarrow \GM^r$ is monomial. 
By elementary linear algebra, there exists a group morphism
$\varrho : \GM^r \rightarrow \GM^r$ such that if we compose the original $\GM^r$-action
on $Y$ with $\varrho$, the morphism $ \pi$ becomes  diagonally monomial
for that new action. 
Furthermore, the 
image
of
$(p, \pi) : Y \rightarrow S \times \GM^r$ with the action twisted by $\varrho$
in 
$\Var_{S \times \GM^r}^{\GM^r}$, hence also its class
in $K_0 (\Var_{S \times \GM^r}^{\GM^r})$
and in $\cM_{S \times \GM^r}^{\GM^r}$,  does not
depend on $\varrho$. We denote that class by
$[(p, \pi): Y \rightarrow S \times \GM^r]$.
Indeed, the first statement  amounts to saying that
for every matrix $A$ in ${\rm M}_r (\ZZ) \cap {\rm GL}_r  (\QQ)$
there exists $B$ in  ${\rm M}_r (\ZZ) \cap {\rm GL}_r (\QQ)$ such that 
$BA$ is diagonal with coefficients in $\NN_{>0}$, and the second one
follows from the observation that if $B'$ is another such matrix,
there exists diagonal matrices $C$ and $C'$ with coefficients in $\NN_{>0}$
such that $CB = C' B'$.

More generally if $W$ is a constructible subset 
of Y stable by the $\GM^r$-action, we shall call a morphism
$\pi : W \rightarrow \GM^r$ piecewise monomial if 
there is a finite
partition of $W$ into locally closed
$\GM^r$-invariant subsets on which the restriction of $\pi$
is a
monomial morphism. 
To such a $W$ endowed with a morphism $(p, \pi) : W \rightarrow S \times \GM^r$
such that the
fibers of $p : W \rightarrow S$ are
$\GM^r$-invariant and $\pi:  W
\rightarrow \GM^r$ is piecewise monomial,  we assign by additivity a
class $[(p, \pi): W \rightarrow S \times \GM^r]$ in $\cM_{S \times \GM^r}^{\GM^r}$.

\subsection{Rational series}\label{sr}Let $A$ be one of the rings
$\ZZ [\LL, \LL^{-1}]$,
$\ZZ [\LL, \LL^{-1}, (\frac{1}{1 - \LL^{-i}})_{i > 0}]$,
$\cM_{S \times \GM^r}^{\GM^r}$.
We denote by
$A[[T]]_{\rm sr}$
the
$A$-submodule of
$A [[T]]$
generated 
by $1$ and by finite
products of terms
$p_{e, i} (T) = \frac{\LL^e T^i}{1 - \LL^e T^i}$,
with $e$ in $\ZZ$ and $i$ in $\NN_{>0}$.
There is a unique
$A$-linear morphism
\begin{equation}
\lim_{T \mapsto \infty} : A[[T]]_{\rm sr} \longrightarrow A
\end{equation}
such that
\begin{equation}
\lim_{T \mapsto \infty} (\prod_{i \in I} p_{e_i, j_i} (T)) = (- 1)^{|I|},
\end{equation}
for every family $((e_i, j_i))_{i \in I}$ in
$\ZZ \times \NN_{>0}$,
with $I$ finite, maybe empty.

\subsection{}\label{chi}Let $I$ be a finite set. We shall consider rational
polyhedral convex cones in $\RR_{>0}^I$.
By this we mean a convex subset of
$\RR_{>0}^I$
defined by a finite number of integral linear inequalities
of type
$a \geq 0$ or $b > 0$, and
stable by multiplication by $\RR_{>0}$.
Let $\Delta$ be such a cone in
$\RR_{>0}^I$. We denote by
$\bar \Delta$ its closure
in $\RR_{\geq 0}^I$.

Let $\ell$ and  $\nu$ be  integral
linear forms on $\ZZ^I$ which are  positive  on $\bar{\Delta}\setminus \{0\}$.
Let us consider the series
\begin{equation}
S_{ \Delta, \ell, \nu} (T):= \sum_{k\in {\Delta} \cap
\NN_{>0}^I} T^{\ell(k)}\LL^{- \nu(k)}
\end{equation}
in  $\ZZ [\LL, \LL^{-1}] [[T]]$.

In the special case when $\Delta$ is open
in its in linear span
and
$\bar \Delta$
is generated by vectors
$(e_1,\ldots,e_m)$ which are part of a $\ZZ$-basis of
the $\ZZ$-module $\ZZ^I$,
the series $S_{ \Delta, \ell, \nu}$
lies in $\ZZ [\LL, \LL^{-1}] [[T]]_{\rm sr}$
and $\lim_{T \mapsto \infty} S_{ \Delta, \ell, \nu} (T)$ is equal to $(-1)^{\dim(\Delta)}$.
By additivity with respect to disjoint union of
cones with the positivity assumption,
one deduces that, in general,
$S_{ \Delta, \ell, \nu}$
lies in $\ZZ [\LL, \LL^{-1}] [[T]]_{\rm sr}$
and $\lim_{T \mapsto \infty} S_{ \Delta, \ell, \nu} (T)$ is equal to $\chi (\Delta)$,
the Euler characteristic with compact supports of $\Delta$.

In particular we get the following lemma
(compare with Lemma 2.1.5 in
\cite{guibert} and \cite{D89} pp.
1006-1007):

\begin{lem}\label{geomsum}Let $\Delta$ be a rational
polyhedral convex cone
in $\RR_{>0}^{I}$ defined by  
\begin{equation}\sum_{i \in K} a_i
x_i \leq \sum_{i \in I \setminus K} a_i x_i,
\end{equation} with $a_i$
in $\NN$ and $a_i >0$, for $i$ in $K$, and with $K$ and $I
\setminus K$ non empty.
If $\ell$ and $\nu$ are integral linear forms positive on
$\bar \Delta \setminus \{0\}$, 
then $\lim_{T \mapsto \infty} S_{ \Delta, \ell, \nu} (T) = 0$.
 \qed
\end{lem}

\section{Motivic vanishing cycles}\label{mvc}

\subsection{Arc spaces} We denote as usual by  $\cL_n (X)$
the space of arcs of order $n$,  also
known as the  $n$-th jet space on $X$.
It is a $k$-scheme whose $K$-points, for $K$ a field containing $k$,
is the set of morphisms
$\varphi : \Spec K [t]/ t^{n + 1} \rightarrow X$.
There are canonical morphisms $\cL_{n + 1} (X) \rightarrow
\cL_n (X)$ which are $\AA^d_k$-bundles when $X$ is smooth of pure dimension $d$.
The
arc space $\cL (X)$ is
defined as the projective limit of this system.
We denote by $\pi_n : \cL (X) \rightarrow \cL_n (X)$
the canonical morphism.
There is a canonical $\GM$-action on $\cL_n (X)$ and on $\cL (X)$
given by
$a \cdot  \varphi (t) = \varphi (at)$.

For an element $\varphi$ in $K [[t]]$ or in $K[t] /t^{n + 1}$, we
denote by $\ord_t (\varphi)$ the valuation of $\varphi$ and
by
$\ac (\varphi)$ its first non zero coefficient, with the convention
$\ac (0) = 0$.

\subsection{Motivic zeta function and Motivic Milnor fiber}\label{mz}
Let us start by recalling some basic
constructions
introduced by
Denef and Loeser in \cite{motivic}, \cite {lef} and \cite{barc}.

Let $X$ be a smooth variety over $k$ of pure dimension $d$ and
$g : X \rightarrow \AA^1_k$.
We set $X_0 (g)$ for the zero locus of $g$,
and consider, for $n \geq 1$,
the variety
\begin{equation}
\cX_n (g) := \Bigl \{\varphi \in \cL_n (X)  \Bigm  \vert \ord_t g (\varphi) = n \Bigr\}.
\end{equation}
Note that $\cX_n (g)$ is invariant by the $\GM$-action on $\cL_n (X)$.
Furthermore $g$ induces a morphism
$g_n  : \cX_n (g) \rightarrow \GM$, assigning to a point
$\varphi$ in $\cL_n (X)$
the coefficient 
$\ac (g (\varphi))$
of $t^n$ in $g (\varphi)$, which we shall also denote by $\ac (g) (\varphi)$.
This morphism is diagonally monomial of weight $n$
with respect to the $\GM$-action
on $\cX_n (g)$ since $g_n (a \cdot \varphi) = a^n g_n (\varphi)$,
so we can consider the class $[\cX_n (g)]$ of
$\cX_n (g)$ in $\cM_{X_0 (g) \times \GM}^{\GM}$.

We now consider the motivic zeta function
\begin{equation}
 Z_g(T): =\sum_{n\geq 1}[\cX_n(g)]\,\mathbf{L}^{-nd}\,T^n
\end{equation}
in $\cM_{X_0 (g) \times \GM}^{\GM}[[T]]$.
Note that
$Z_g = 0$  if $g = 0$ on $X$.

Denef and Loeser showed in \cite{motivic} and \cite{barc} (see also \cite{lef})
that $Z_g(T)$ is a rational series by giving a formula for
$Z_g(T)$ in terms of a resolution of $f$.

\subsection{Resolutions}\label{resolutions}
Let us introduce some notation and terminology. Let $X$ be
a smooth variety of pure dimension $d$ and let $Z$ a closed
subset of $X$ of codimension everywhere $\geq 1$. By a
log-resolution $h : Y \rightarrow X$ of $(X, Z)$, we mean a
proper morphism $h : Y \rightarrow X$ with $Y$ smooth such
that the restriction of $h : Y \setminus h^{-1} (Z)
\rightarrow X \setminus Z$ is an isomorphism, and $h^{-1}
(Z)$ is a divisor with normal crossings. We denote by
$E_i$, $i$ in $A$, the set of irreducible components of the
divisor $h^{-1}  (Z)$. For $I \subset A$, we set 
\begin{equation}E_I := \bigcap_{i \in I} E_i
\end{equation}and
\begin{equation}
E_I^{\circ} := E_I \setminus \bigcup_{j \notin I} E_j.
\end{equation}
We denote by $\nu_{E_i}$  the normal bundle of $E_i$ in
$Y$, by $\nu_{E_I}^J$ the fiber product, for $J$ contained
in $I$, of the restrictions to $E_I$ of the bundles
$\nu_{E_i}$, $i$ in $J$, and by $\pi_I^J :
\nu_{E_I}^J\rightarrow E_I$ the canonical projections. For
any of these vector bundles $\nu$ we will denote by
$\overline{\nu}$ the projective bundle associated to the
sum of $\nu$ with the trivial line bundle.

We will denote by $U_{E_i}$ the complement of the zero
section in $\nu_{E_i}$ and by $U_I^J$ (resp. $U_{E_I}^J$) the fiber product,
for $J$ contained in $I$, of the restrictions of the spaces
$U_{E_i}$, $i$ in $J$, to $E_I^{\circ}$ (resp. $E_I$). We will still
denote by $\pi_I^J$ the induced projection from $U_I^J$ (resp. $U_{E_I}^J$) 
onto $E_I^\circ$ (resp. $E_I$).

When $J=I$ we will simply write  $\nu_{E_I}$ (resp.
$\overline{\nu}_{E_I}$, $\pi_I$, $U_I$, $U_{E_I}$) for
$\nu_{E_I}^I$ (resp. $\overline{\nu}_{E_I}^I$, $\pi_I^I$,
$U_I^I$, $U_{E_I}^I$).

\bigskip
If $\cI$ is a sheaf of ideals defining a closed subscheme $Z$ 
and $h^{-1} (\cI)\cO_Y$ is locally principal,
we define $N_i (\cI)$, the multiplicity
of $\cI$ along $E_i$,
by the equality of divisors
\begin{equation}
h^{-1} (Z) = \sum_{i \in A} N_i (\cI) E_i.
\end{equation}
If $\cI$ is principal generated by a function $g$ we write
$N_i (g)$ for $N_i (\cI)$.
Similarly,  we define integers $\nu_i$ by
the equality of divisors
\begin{equation}
K_Y = h^* {K_X} + \sum_{i \in A} (\nu_i - 1) E_i.
\end{equation}

Let $\cI_1$ and $\cI_2$ be two sheaves of ideals on $X$
whose associated reduced closed subschemes $Z_1$ and $Z_2$
have codimension at least one.
Let $h : Y \rightarrow X$
be a log-resolution of
$(X, Z_1 \cup Z_2)$ such that $h^* (\cI_1)$ and
$h^* (\cI_2)$ are locally principal.
Then we set
\begin{equation}\label{teq}
 \gamma_h (\cI_1, \cI_2)  := \sup_{\{i \in A |N_i (\cI_2) >0\}} \frac{N_i (\cI_1)}{N_i (\cI_2)}.
\end{equation}
If $x$ is a closed point of $Z_2$, 
we set
\begin{equation}\label{teqvar}
\gamma_{h, x} (\cI_1, \cI_2)  = \sup_{\{i \in A_x |N_i (\cI_2)
>0\}}  \frac{N_i (\cI_1)}{N_i (\cI_2)},
\end{equation}
 with $A_x$ the set of $i$ in $A$ such that
$\vert h^{-1}Ê(x)\vert \cap E_i \not= \emptyset$.
Finally we define 
$\gamma (\cI_1, \cI_2)$, resp.
$\gamma_x (\cI_1, \cI_2)$
as the infimum of all 
$\gamma_h (\cI_1, \cI_2)$, resp.  $\gamma_{h, x} (\cI_1, \cI_2)$,
for $h$ a log-resolution of
$(X, Z_1 \cup Z_2)$ such that $h^* (\cI_1)$ and
$h^* (\cI_2)$ are locally principal.

\subsection{}\label{g_I}
Let $g$ be a function on a smooth variety $X$ of
pure dimension $d$. 
Assume $X_0 (g)$ is nowhere dense
in $X$.
Let $F$ a reduced divisor containing
$X_0 (g)$ and let $h : Y \rightarrow X$ be a log-resolution
of $(X, F)$. 
We fix $I$ such that
there exists $i$ in $I$ with $N_i (g) > 0$. 
Let us explain how $g$ induces a morphism $g_I
: U_I \rightarrow \GM$. Note that the function $g \circ h$
induces a function
\begin{equation}\label{gtt}
\bigotimes_{i \in I} \nu_{E_i}^{\otimes N_i (g)}{}_{|E_I}
\longrightarrow \AA^1_k.
\end{equation}
We define $g_I
: \nu_{E_I} \rightarrow \AA^1_k$ as the composition of this
last function with the natural morphism $\nu_{E_I}
\rightarrow \otimes_{i \in I} \nu_{E_i}^{\otimes N_i
(g)}{}_{|E_I}$, sending
$(y_i)$ to $\otimes y_i^{\otimes
N_i (g)}$. 
We still denote by $g_I$ the induced morphism
from $U_I$ (resp. $U_{E_I}$) to $\GM$ (resp. $\AA^1_k$).

We view $U_I$ as a variety over $X_0 (g) \times \GM$ via
the morphism $(h\circ\pi_I,g_I)$. The group $\GM$ has a
natural action on each $U_{E_i}$, so the diagonal action
induces a $\GM$-action on $U_I$. Furthermore, the morphism
$g_I$ is
monomial,
so $ U_I \rightarrow X_0 (g) \times
\GM$ has a class in $\cM_{X_0 (g) \times \GM}^{\GM}$ which
we will denote by $[U_I]$.

\subsection{}\label{desdef}The morphism $g_I$ may be described in terms of 
the following variant of the deformation to the normal cone to $E_{I}$
in
$Y$, 
cf. \cite{ful}.  We consider
the affine space $\AA^I_k = \Spec k [u_i]_{ i\in I}$ and the
subsheaf
\begin{equation}\label{grd}
\cA_{I}:= \sum_{{\bf n} \in \NN^I} \cO_{Y \times \AA^I_k}
\left(- \sum_{i \in I} n_i (E_i \times \AA^I_k) \right )
\prod_{i \in I}u_i^{-n_i}
\end{equation}
of $\cO_{Y \times \AA^I_k} [u_i^{-1}]_{ i\in I}$. It is a
sheaf of rings and we set $
CY_{I} := \Spec \cA_{I}.$ 
The natural inclusion $\cO_{Y \times \AA^I_k} \rightarrow
\cA_{I}$ induces a morphism $\pi : CY_{I} \rightarrow Y
\times \AA^I_k$, hence a morphism $p : CY_{I} \rightarrow
\AA^I_k$. 
The ring $\cA_I$ being a graded subring of the ring $\cO_Y
[u_i, u_i^{-1}]_{ i\in I}$, we consider the corresponding
$\GM^I$-action $\sigma_I$ on $CY_{I}$, leaving sections of
$\cO_Y$ invariant and acting by $(\lambda_i, u_i)\mapsto
\lambda_i^{-1} u_i$ on $u_i$. 
We may then identify equivariantly $\nu_{E_I}$ with  the fiber $p^{-1}Ê(0)$.
The image by the inclusion $\cO_{Y \times \AA^I_k} \rightarrow
\cA_{I}$ of the function $g\circ h$
is divisible by $ \prod_{i \in I}u_i^{N_i(g)}$ in
$\cA_{I}$, so we may consider the quotient $\widetilde g_{I}$ in
$\cA_{I}$.
The restriction of  $\widetilde g_{I}$ to the fiber
$p^{-1}(0) \simeq \nu_{E_I}$ is nothing else than $g_I$.
As $g$ may vanish only on the divisors $E_i$, $i$ in $A$, 
the function $g_I$ does not vanish on $U_I$ and induces
a monomial morphism $g_I: U_I \longrightarrow \GM$.

Let us note the following ``transitivity'' property.
If we write $I$ as  a disjoint union $K\sqcup J$, one notices that 
$p^{-1}(0\times \GM^J)$ is equivariantly isomorphic to  
$\nu_{E_K}\times \GM^J$. Hence, restricting
$p : CY_I \rightarrow \AA_k^I$ to 
$p^{- 1}Ê(0 \times \AA_k^J)$, 
the function $g_I: U_I \longmapsto \GM$ can be obtained from $g_K$
by the same process as we obtained it from $g$,
replacing $Y$ by $\nu_{E_K}$, $I$ by $J$ and
$g$ by $g_K:\nu_{E_K}\longrightarrow \AA^1_k$.

\subsection{}We now assume that $F = X_0 (g)$, that
is $h : Y \rightarrow X$ is a log-resolution of $(X,
X_0(g))$. In this case, $h$ induces a bijection between
$\cL (Y) \setminus \cL (\vert h^{-1} (X_0 (g))\vert)$ and
$\cL (X) \setminus \cL (X_0 (g))$. One deduces by using the
change of variable formula, in a way completely similar to
\cite{motivic} and \cite{barc}, the equality
\begin{equation}\label{dl1}
Z_g (T)  = \sum_{\emptyset \not= I \subset A} [U_I]
\prod_{i \in I} \frac{1}{T^{- N_i (g)} \LL^{\nu_i} - 1}
\end{equation}
in $\cM_{X_0 (g) \times \GM}^{\GM}[[T]]$.

In particular, the function  $Z_g (T)$ is rational
and belongs to
$\cM_{X_0 (g) \times \GM}^{\GM} [[T]]_{\rm sr}$,
with the notation of \ref{sr}, hence we can consider
$ \lim_{T \mapsto \infty} Z_g (T)$
in
$\cM_{X_0 (g) \times \GM}^{\GM}$.

We set
\begin{equation}
\cS_g := - \lim_{T \mapsto \infty} Z_g (T),
\end{equation}
which by (\ref{dl1})  may be expressed on a resolution $h$ as
\begin{equation}\label{dl2}
\cS_g  = - \sum_{\emptyset \not=  I \subset A} (-1)^{|I|}[U_I],
\end{equation}
in $\cM_{X_0 (g) \times \GM}^{\GM}$.

We shall also consider in this paper
the motivic
vanishing cycles defined
as
\begin{equation}
\cS_g^{\phi} := (-1)^{d-1}(\cS_g - [\GM \times X_0 (g)])
\end{equation}
in $\cM_{X_0 (g) \times \GM}^{\GM}$.
Here $d$ denotes the dimension of
$X$ and $\GM \times X_0 (g)$
is endowed with the
standard $\GM$-action on the first factor and the
trivial $\GM$-action on the second factor.

\subsection{A modified zeta function}
We now explain how to extend $\cS_g$ to the whole
Grothendieck group $\cM_{X}$ in such a way that $\cS_g ([X
\rightarrow X])$ is equal to $\cS_g$. A similar result has been obtained
by F. Bittner in \cite{bittner2}. We present here a
somewhat different approach that avoids the use of the weak
factorization Theorem, by constructing directly $\cS_g ([Y
\rightarrow X])$ for generators of $\cM_{X}$ of the form $Y
\rightarrow X$ with $Y$ smooth.

Let $X$ be a smooth variety of pure dimension $d$ and
let $U$ be a dense open in $X$.
Consider again a function
$g : X \rightarrow \AA^1_k$.
We denote by $F$ the closed subset
$X \setminus U$ and by $\cI_F$
the ideal of functions vanishing on $F$.
We start by defining
$\cS_g( [U \rightarrow X])$.

Fix $\gamma \geq 1$ a positive integer.
We will consider the modified zeta function
$Z^{\gamma}_{g, U} (T)$
defined as follows.
For $n \geq 1$, we consider
the constructible set
\begin{equation}\label{grm}
\cX_n^{\gamma n} (g, U) := \Bigl \{\varphi \in \cL_{\gamma
n} (X)  \Bigm  \vert \ord_t g(\varphi) = n, \ord_t
\varphi^*(\cI_F) \leq \gamma n \Bigr\}.
\end{equation}
As in \ref{mz}, we consider the morphism
$\cX_n^{\gamma n} (g, U) \rightarrow \GG_m$ induced by
$\varphi \mapsto \ac (g (\varphi))$.
It is piecewise monomial, so we can consider
the 
class
$[\cX^{\gamma n}_n (g, U)]$ in
$\cM_{X_0 (g) \times \GM}^{\GM}$ by \ref{twist}.
We set
\begin{equation}\label{mzf}
Z^{\gamma}_{g, U} (T) := \sum_{n\geq 1}[\cX^{\gamma n}_n
(g, U)]\,\mathbf{L}^{-  \gamma  n d}\,T^n
\end{equation}
in $\cM_{X_0 (g) \times \GM}^{\GM} [[T]]$. Note that for $U = X$, $Z^{\gamma}_{g, U} (T)$ is equal to $Z_g (T)$ for
every $\gamma$, since in this case, $[\cX^{\gamma n}_n (g,
U)]  \LL^{- \gamma n d}= [\cX_n (g)] \LL^{-n d}$. 
Note also that $Z^{\gamma}_{g, U} (T) = 0$ if
$g$ is identically zero on $X$.

If $X_0 (g)$ is nowhere dense in $X$ and
$h : Y \rightarrow X$ is a log-resolution of $(X, F\cup
X_0 (g))$, we denote by $C$ the set $\{ i \in A \mid N_i
(g) \not=0\}$.

\begin{prop}\label{formop}Let $U$ be a dense open in
the smooth variety $X$ of pure dimension $d$ with a function
$g : X \rightarrow \AA^1_k$.
There exists $\gamma_0$ such that for every $\gamma > \gamma_0$
the series $Z^{\gamma}_{g, U} (T)$ lies
in
$\cM_{X_0 (g) \times \GM}^{\GM} [[T]]_{\rm sr}$
and
$\lim_{T \mapsto \infty} Z^{\gamma}_{g, U} (T)$
is independent of
$\gamma > \gamma_0$.
We set
$\cS_{g, U} = -\lim_{T \mapsto \infty} Z^{\gamma}_{g, U} (T)$.
Furthermore, 
if $X_0 (g)$ is nowhere dense in  $X$ and 
$h : Y \rightarrow X$ is a log-resolution
of
$(X, F \cup X_0 (g))$,
\begin{equation}\label{fty}
\cS_{g, U} = - %
\sum_{\sur{I \not= \emptyset}{I \subset C}}
(-1)^{\vert I \vert}\,[U_I] = h_!\left(\cS_{g\circ
h,h^{-1}(U)}\right)
\end{equation}
in $\cM_{X_0 (g) \times \GM}^{\GM}$.
\end{prop}

\begin{proof}We may assume $X_0 (g)$ is nowhere dense in $X$.
Let $h : Y \rightarrow X$ be a log-resolution
of
$(X, F \cup X_0 (g))$.
As in  the proof of Theorem 2.4 of \cite{lef},
we deduce from the change of variable formula,
or more precisely from Lemma 3.4 in \cite{arcs},  that
\begin{equation}\label{trtr}Z^{\gamma}_{g, U} (T)  = \sum_{I \cap C \not= \emptyset}
[U_I] \, S_I (T)
\end{equation}
with
\begin{equation}
S_I (T) = \sum_{\sur{k_i \geq 1}{
\sum_i k_i N_i (\cI_F) \leq \gamma \sum_i k_i N_i (g)}}
\prod_{i \in I} (T^{N_i (g)} \LL^{- \nu_i})^{k_i}.
\end{equation}
Assume first that $I \subset C$.
For $\gamma \geq \sup_{i \in I}\frac{N_i (\cI_F)}{N_i (g)}$,
we have  $\sum_i k_i N_i (\cI_F) \leq \gamma \sum_i k_i
N_i (g)$ for all  $k_i\geq 1$, $i\in I$.
It follows that $S_I (T)$ 
lies in $\cM_{X_0 (g) \times \GM}^{\GM} [[T]]_{\rm sr}$
and
$\lim_{T \mapsto \infty} S_I (T) = (-1)^{|I|}$, as soon
as $\gamma \geq \sup_{i \in I}\frac{N_i (\cI_F)}{N_i (g)}$.

Now assume
$ \emptyset\neq I\setminus C =K$.
For $\gamma \geq \sup_{i \in I\setminus K}\frac{N_i (\cI_F)}{N_i
(g)}$,
the sum runs over the points with
coordinates in $\NN_{>0}$ of the cone $\Delta_I$ 
in
$\RR_{>0}^I$
defined by the single
inequality
\begin{equation}
\sum_{i \in K} a_i x_i \leq \sum_{i \in I \setminus K} a_i x_i,
\end{equation}
with
$a_i$ in $\NN$ and $a_i >0$,
for $i$ in $K$. 
Note that both $K$ and $I \setminus K$ are non empty.
It follows from Lemma \ref{geomsum}
that in this case $S_I (T)$ 
lies in $\cM_{X_0 (g) \times \GM}^{\GM} [[T]]_{\rm sr}$
and
$\lim_{T \mapsto \infty} S_I (T) = 0$.
The statement we have to prove then holds if
we set
$\gamma_0 =  \sup_{i \in  C} 
\frac{N_i (\cI_F)}{N_i (g)} 
= \gamma_h (\cI_F, (g))$.
Note that since this holds for any $h$, we could also take 
$\gamma_0 = \gamma (\cI_F, (g))$.
\end{proof}

\begin{theorem}[Extension to the Grothendieck group]\label{ext}Let $X$ be a variety
with a function
$g : X \rightarrow \AA^1_k$.
There exists a unique $\cM_k$-linear group morphism
\begin{equation}\label{exteq}
\cS_g : \cM_X \longrightarrow \cM_{X_0 (g) \times \GM}^{\GM}
\end{equation}
such that, for every proper morphism $p : Z \rightarrow X$,
with $Z$ smooth, and every dense open subset $U$ in $Z$,
\begin{equation}\label{skk}
\cS_g ([U \rightarrow X]) = p_! (\cS_{g \circ p, U}).
\end{equation}
\end{theorem}

\begin{proof}Since $K_0 (\Var_X)$ is generated by classes
$[U \rightarrow X]$ with $U$ smooth 
connected
and every such  $U
\rightarrow X$ may be embedded in a proper morphism $Z
\rightarrow X$ with $Z$ smooth and $U$ dense in $Z$,
uniqueness is clear. For existence let us first note that
if we define $\cS_g ([U \rightarrow X]) = \cS_g ([U])$ by
the  right hand side of (\ref{skk}), the result is
independent from the choice of the embedding  in a proper
morphism $p : Z \rightarrow X$. Indeed, 
this is clear if $g \circ p$ vanishes identically on $U$,
so we may assume $(g \circ p)^{-1} ( 0)$ is of codimension $>0$.
In this case,
if we have another
such morphism $p' : Z' \rightarrow X$, there exists a
smooth variety $W$ with proper morphisms $h : W \rightarrow
Z$ and $h' : W \rightarrow Z'$, such that $p \circ h = p'\circ h'$ and $h$ and $h'$ are respectively
log-resolutions of $(Z, (Z \setminus U) \cup (g \circ
p)^{-1} (0))$ and $(Z', (Z' \setminus U) \cup (g \circ
p')^{-1} (0))$, so the statement follows from (\ref{fty}).

Let us now prove the following additivity
statement: if $\kappa :  U \rightarrow X$ is a morphism with
$U$ smooth and $W$ is a smooth closed subset of $U$,
then
\begin{equation}\label{oiu}
\cS_g ([U \rightarrow X]) = \cS_g ([W \rightarrow X]) +
\cS_g ([U \setminus W \rightarrow X]).
\end{equation}
We may assume $U$ and $W$ are connected and $U \setminus W$ is dense in $U$.
The result being trivial if $g \circ \kappa$ vanishes identically, we may assume
this is not the case.
By Hironaka's
strong resolution  of singularities, we may embed
$U$ in a smooth variety
$Z$ with $p : Z \rightarrow X$
a proper morphism extending $\kappa$ such that
$Z \setminus U$ is a normal
crossings divisor and
the closure  $\overline W$ of $W$
in $Z$ is smooth. 
Again by Hironaka's
strong resolution  of singularities,
there exists a log-resolution
$h : \widetilde Z \rightarrow Z$
of $(Z, (Z \setminus U) \cup (g \circ p)^{-1}Ê(0))$ such that
the closure $\widetilde W$ of $h^{-1}Ê(W)$ in
$\widetilde Z$ is smooth and intersects the divisor
$D := h^{-1} ((Z \setminus U) \cup (g \circ p)^{-1}Ê(0))$ transversally.
We denote by
$E_i$, $i$ in $A$, the irreducible components of the
divisor $D$ and use the notations of \ref{resolutions}.
It follows from the definition and (\ref{fty})
that
\begin{equation}\label{one}
\cS_{g} ([U \rightarrow X]) =  - %
 \sum_{\sur{I \not= \emptyset}{I \subset C}}
(-1)^{\vert I \vert}\,[U_I],
\end{equation}
in $\cM_{X_0 (g) \times \GM}^{\GM}$.
Note that if $W$ is contained in 
$(g \circ p)^{-1}Ê(0)$, the above discussion still holds
for $U$ replaced by $U \setminus W$, so we have
$\cS_{g} ([U \rightarrow X]) = \cS_{g} ([(U \setminus W) \rightarrow X])$,
and  (\ref{oiu}) 
follows, 
since $\cS_{g} ([W \rightarrow X]) = 0$ in this case.

Now we assume $W$ is not contained in 
$(g \circ p)^{-1}Ê(0)$.
Note 
that 
the morphism
$h_0 : \widetilde W \rightarrow \overline W$ induced by $h$ is a log-resolution
of $(\overline W, (\overline W \setminus W) \cup (g \circ p)_{\vert \overline W}^{-1}Ê(0))$.
Furthermore the irreducible components of the 
normal crossings
divisor $h_0^{-1} ((\overline W \setminus W) \cup (g \circ p)_{\vert \overline W}^{-1}Ê(0))$
are exactly those amongst the $E_i \cap \widetilde W$ which are non empty.
Hence, denoting by
$U_I \vert_{E^{\circ}_I \cap \widetilde W}$
the restriction of the bundle $U_I$ to $E^{\circ}_I \cap \widetilde W$,
it follows from the definition and (\ref{fty})
that
\begin{equation}\label{two}
\cS_{g} ([W \rightarrow X]) = - %
\sum_{\sur{I \not= \emptyset}{I \subset C}}
(-1)^{\vert I \vert}\,[U_I \vert_{E^{\circ}_I \cap \widetilde W}],
\end{equation}
in $\cM_{X_0 (g) \times \GM}^{\GM}$.
Let us now consider
the blowing up $h' : Z' \rightarrow \widetilde Z$ of $\widetilde Z$ along $\widetilde W$.
The 
exceptional divisor
 $W'$ of $\widetilde W$ is smooth. Furthermore
$h \circ h' : Z' \rightarrow Z$
is a log-resolution of $(Z, (Z \setminus (U \setminus W)) \cup (g \circ p)^{-1}Ê(0))$,
and $D' := (h \circ h')^{-1} ((Z \setminus (U \setminus W)) \cup (g \circ p)^{-1}Ê(0))$
is a normal crossings divisor whose irreducible components
are the strict transforms $E'_i$ of $E_i$ in $Z'$, $i$ in $A$ together with
$W'$. We set $A' := A \sqcup \{0\}$ and $E'_0 := W'$ in order 
to use the notations of \ref{resolutions} in this setting, adding everywhere  $'$ as an exponent.
Again, it follows from the definition and (\ref{fty})
that
\begin{equation}\label{three}
\cS_{g} ([(U \setminus W) \rightarrow X] ) = - %
\sum_{\sur{I \not= \emptyset}{I \subset C'}}
(-1)^{\vert I \vert}\,[U'_I],
\end{equation}
in $\cM_{X_0 (g) \times \GM}^{\GM}$,
where $C' = \{i \in A' \vert N_i (g \circ p \circ h \circ h') \not= 0\}$.
The hypothesis made on $W$ insures that $C' = C$.
So it 
is
enough to prove that for $I$ non empty and contained in $C$,
\begin{equation}\label{four}
[U_I ] = [U_I \vert_{E^{\circ}_I \cap \widetilde W}] +
[U'_I]
\end{equation}
in $\cM_{X_0 (g) \times \GM}^{\GM}$, 
which follows from the fact that
the restriction 
$U_I \vert_{E^{\circ}_I \setminus (E^{\circ}_I \cap \widetilde W)}$
of the bundle $U_I$ to $E^{\circ}_I \setminus (E^{\circ}_I \cap \widetilde W)$
and the bundle $U'_I$ have the same class in 
$\cM_{X_0 (g) \times \GM}^{\GM}$, since $h'$ is an isomorphism outside
$W'$. This concludes the proof of (\ref{oiu}).

Let again $U \rightarrow X$ be in $\Var_X$ with $U$ smooth and connected.
Let $W$ be a smooth proper variety over $k$.
Note that 
\begin{equation}\label{monk}
\cS_g ([W \times U \rightarrow X ]) = [W] \, \cS_g ([U \rightarrow X])
\end{equation}
in $\cM_{X_0 (g) \times \GM}^{\GM}$. 
Indeed, let us embed $U \rightarrow X $ in $p : Z \rightarrow X$
with $Z$ smooth and proper and $U$ dense in $Z$.
We may assume $g \circ p$ is not identically zero.
If
$h : Y \rightarrow Z$ is a log-resolution
of $(Z, (Z \setminus U) \cup (g \circ p)^{-1} (0))$,
then $W \times U \rightarrow X $ may be embedded  in $W \times Z \rightarrow X$,
 ${\rm id}Ê\times h : W \times Y \rightarrow W \times Z$
is a log-resolution of $(W \times Z, ((W \times Z) \setminus (W \times U)) \cup (W \times
g \circ p)^{-1} (0))$, hence (\ref{monk})
follows from (\ref{fty}) and (\ref{skk}).
By the additivity statement we already proved, relation (\ref{monk}) in fact holds for
every variety $W$ over $k$, so our construction of 
$\cS_g$ may be extended uniquely by $\cM_k$-linearity to
a  $\cM_k$-linear group morphism
$\cM_X \rightarrow \cM_{X_0 (g) \times \GM}^{\GM}$, which finishes the proof.
\end{proof}

\subsection{The equivariant setting}\label{eqset}
Let $X$ be a variety with a function $g : X \rightarrow
\AA^1_k$. By Theorem \ref{ext}, there is a canonical
morphism
\begin{equation}
\cS_g : \cM_{X} \longrightarrow \cM_{X_0 (g)\times \GM}^{\GM}.
\end{equation}

We want to lift this morphism
to a morphism, still denoted by $\cS_g$,
\begin{equation}\label{equv}
\cS_g : \cM_{X \times \GM^r}^{\GM^r} \longrightarrow
\cM_{X_0 (g) \times \GM^r \times \GM}^{\GM^r  \times \GM}
\end{equation}
such that the diagram
\begin{equation}\label{no}\xymatrix{
\cM_{X \times \GM^r}^{\GM^r}
\ar[d] \ar[r]^{\cS_g} &
\cM_{X_0 (g) \times \GM^r \times \GM}^{\GM^r  \times \GM}
 \ar[d]\\
\cM_X
\ar[r]^{\cS_g}&\cM_{X_0 (g)\times \GM}^{\GM}
}
\end{equation}
is commutative, the vertical arrows being
given by forgetting the $\GM^r$-action and taking the fiber over $1$ in $\GM^r$.

Let us start with some basic facts we shall use
without further mention.
We fix the variety $X$ which we shall consider as endowed with the trivial
$\GM^r$-action. Let $Z$ be a smooth
variety of pure dimension $d$ endowed with a good $\GM^r$-action
and an equivariant morphism $p: Z \rightarrow X$. 
The induced action
on the affine bundles $\cL_{n + 1} (Z) \rightarrow \cL_n (Z)$
is affine. In particular, by relation (\ref{addrel}), 
$[\cL_{n + 1} (Z) \rightarrow X] = \LL^d [\cL_n (Z) \rightarrow X]$ in $\cM_{X \times \GM^r}^{\GM^r}$.
Similarly, if $h : Y \rightarrow Z$ is a proper birational
$\GM^r$-equivariant morphism with 
$Y$
smooth with a 
good $\GM^r$-action, the fibrations occuring in Lemma 3.4 of \cite{arcs}
are (piecewise) affine bundles  and the induced 
$\GM^r$-action is affine, hence, by relation (\ref{addrel}),
one does not see the action on the fibers in the Grothendieck ring
 $\cM_{X \times \GM^r}^{\GM^r}$.

We now assume $X$
is endowed with a morphism
$g:X\rightarrow \AA_k^1$
and that $Z$ is endowed with a monomial morphism
$\mathbf{f} = (f_1, \dots, f_r) : Z \rightarrow \GM^r$
such that $(p, \mathbf{f}) : Z \rightarrow X \times \GM^r$ is proper.
We consider an
open
dense subset $U$ of $Z$ stable under the $\GM^r$-action.
Similarly as in (\ref{grm}), we set
\begin{equation}\label{grmbis}
\cX_n^{\gamma n} (g\circ p, U) := \Bigl \{\varphi \in \cL_{\gamma
n} (Z)  \Bigm  \vert \ord_t (g \circ p)(\varphi) = n, \ord_t
\varphi^*(\cI_F) \leq \gamma n \Bigr\},
\end{equation}
with $F := Z \setminus U$.
The $\GM^r$-action on $Z$
induces a $\GM^r$-action on
$\cX^{\gamma n}_n (g \circ p, U)$ via its induced action on the arc space.
On the other side, the standard $\GM$-action on arcs  
\begin{equation}\label{}
(\lambda \cdot \varphi)(t) := \varphi(\lambda t)
\end{equation}
induces a $\GM$-action on
$\cX^{\gamma n}_n (g\circ p, U)$. In this way we get
a $\GM^r \times \GM$-action on $\cX^{\gamma n}_n (g \circ p, U)$.
The morphism 
\begin{equation}(\mathbf{f} \circ \pi_0,\ac(g \circ p))
: \cX^{\gamma n}_n (g \circ p, U)\rightarrow \GM^r\times \GM
\end{equation}
is piecewise monomial,
hence, proceeding as in
\ref{twist}, we may assign to 
\begin{equation}\cX^{\gamma n}_n (g \circ p, U)\longrightarrow X_0(g) \times \GM^r\times \GM
\end{equation}
a class $[\cX^{\gamma n}_n (g\circ p, U)]$
in $\cM_{X_0(g)\times \GM^r\times \GM}^{\GM^r\times \GM}$.

Similarly as in  (\ref{mzf}),
we consider  the corresponding series
\begin{equation}\label{mzfequi}Z^{\gamma}_{g \circ p, U} (T) := \sum_{n\geq 1} \, [\cX^{\gamma n}_n (g \circ p, U)]\,\mathbf{L}^{-
\gamma  n d}\,T^n
\end{equation} in $\cM_{X_0 (g)\times \GM^r \times \GM}^{\GM^r\times \GM}
[[T]]$.

Proceeding  as in the proof of Proposition
\ref{formop}, one proves  that there exists a $\gamma_0$ such that for every
$\gamma > \gamma_0$ the series $Z^{\gamma}_{g \circ p, U} (T)$ belongs to
$\cM_{X_0 (g)\times \GM^r \times \GM}^{\GM^r\times \GM}
[[T]]_{\rm sr}$ and that $\lim_{T \mapsto \infty}
Z^{\gamma}_{g \circ p, U} (T)$ is independent of $\gamma >
\gamma_0$. Indeed, 
we may assume that the zero locus  $Z_0 (g\circ p)$ of $g \circ p$ is nowhere dense in $Z$ and in this case we
now use
a $\GM^r$-equivariant log-resolution of $(Z,
(Z\setminus U) \cup Z_0 (g\circ p))$. (For the existence of
equivariant  resolutions, see \cite{bm} \cite{eh} \cite{ev}
\cite{v1} \cite{v2}.) We now define
$\cS_{g\circ p, U}$ in $\cM_{X_0 (g)\times \GM^r \times \GM}^{\GM^r \times \GM}$
as  $- \lim_{T \mapsto \infty}
Z^{\gamma}_{g\circ p, U} (T)$ for $\gamma >
\gamma_0$.

Still assuming $Z_0 (g\circ p)$ is nowhere dense in $Z$,
let  $h : Y \rightarrow Z$ be such a $\GM^r$-equivariant
log-resolution.
We shall  use  again the notations introduced in \ref{resolutions}.
By connectedness of $\GM^r$,
the $\GM^r$-action on $Z$ induces the trivial action on the set of strata
$E_I^{\circ}$, for $I$ subset
of $A$. 
The $\GM^r$-action on $Y$ induces an action
on the normal  bundles to the divisors $E_{i}$, for
$i$ in $A$, hence on $U_{I}$, for $I$ subset
of $A$.
We also consider the $\GM$-action on $U_{I}$ which is the diagonal action induced by the canonical $\GM^I$-action on $U_I$. In this way we get a $\GM^r\times \GM$-action on $U_{I}$.
Furthermore, with the notation of \ref{g_I},
the morphisms ${\mathbf f}$ and $g$ induce  morphisms
$\mathbf{f}_{I} : U_{I}\rightarrow \GM^r$
and
$\mathbf{g}_{I} : U_{I}\rightarrow \GM$.
Note
that the morphism
$(\mathbf{f}_{I}, g_{I}) : U_{I}\rightarrow \GM^r\times \GM$
is monomial with respect to the $\GM^r\times \GM$-action, since
$g_{I}$ is invariant by the $\GM^r$-action
and monomial with respect to the
$\GM$-action
and the morphism
$\mathbf{f}_{I} : U_{I}\rightarrow \GM^r$
is induced from $\mathbf{f}$ via the projection $U_I \rightarrow Z$.
We can then consider the class  $[U_{I}]$  in
$\cM_{X_0 (g)\times \GM^r \times \GM}^{\GM^r \times \GM}$
of the morphism
\begin{equation}\label{}
(p \circ h\circ\pi_{I}, \mathbf f_{I},g_{I}): U_{I} \longrightarrow
X_0(g)\times\GM^r\times \GM .
\end{equation}
Similarly as in Proposition
\ref{formop}, we get that the equality
\begin{equation}\label{ftyeq}
\cS_{g\circ p, U} = \sum_{I \not= \emptyset, I \subset C}
(-1)^{\vert I \vert}[U_{I}] 
\end{equation}
holds in $\cM_{X_0
(g)\times \GM^r \times \GM}^{\GM^r \times \GM}$.


\begin{remark}When $r = 0$, what is  denoted here by
$[\cX^{\gamma n}_n (g\circ p, U)]$, 
$Z^{\gamma}_{g \circ p, U} (T)$ and
$\cS_{g\circ p, U}$ corresponds to what was denoted by
$p_! ([\cX^{\gamma n}_n (g\circ p, U)])$, 
$p_! (Z^{\gamma}_{g \circ p, U} (T))$ and
$p_! (\cS_{g\circ p, U})$ in the non equivariant setting. This slight conflict of notation 
should lead to  no confusion.
\end{remark}

We can now state the following equivariant analogue of Theorem \ref{ext}.

\begin{theorem}\label{extequi}Let $X$ be a variety
with a function
$g : X \rightarrow \AA^1_k$. 
We consider $X$ endowed with the trivial
$\GM^r$-action.
There exists a unique $\cM_k$-linear group morphism
\begin{equation}\label{equuuuuuu}
\cS_g : \cM_{X \times \GM^r}^{\GM^r} \longrightarrow
\cM_{X_0 (g) \times \GM^r \times \GM}^{\GM^r  \times \GM}
\end{equation}
such that, for every 
smooth variety $Z$ with good $\GM^r$-action,
endowed with an equivariant  morphism
$p  : Z \rightarrow X$
and a monomial morphism ${\mathbf f} : Z \rightarrow \GM^r$ such that the morphism
$(p,  {\mathbf f}) : Z \rightarrow X \times \GM^r$ is proper,
and for every  open
dense subset $U$ of $Z$ which is
stable under the $\GM^r$-action
\begin{equation}\label{skkkkk}
\cS_g ([U \rightarrow X \times \GM^r]) =  
\cS_{g \circ p, U}
\end{equation}
in $\cM_{X_0 (g) \times \GM^r \times \GM}^{\GM^r  \times \GM}$.
\end{theorem}

\begin{proof}Let us denote by  $K'_0 (\Var_{X \times \GM^r}^{\GM^r})$
the Grothendieck ring defined similarly as $K_0 (\Var_{X \times \GM^r}^{\GM^r})$,
but without relation (\ref{addrel}).
Let $U$ be a smooth variety over $k$ with a good
$\GM^r$-action endowed with an equivariant  morphism
$\kappa : U \rightarrow X$, and with a monomial
morphism $\mathbf{f}_U : U \rightarrow \GM^r$.
Note that $U$ may be embedded equivariantly as an open dense
subset of a 
smooth variety $Z$ with good $\GM^r$-action,
endowed with an equivariant  morphism
$p : Z \rightarrow X$ extending $\kappa$ and
a monomial morphism ${\mathbf f} : Z \rightarrow \GM^r$
extending ${\mathbf f}_U$, such that
$(p , {\mathbf f}) : Z \rightarrow X \times \GM^r$ is proper.
Indeed, using the equivalence of categories
of Proposition \ref{equiv} and \ref{twist}, it is enough to know that every smooth variety
$U_0$ endowed with a good $\hat \mu^r$-action and with an equivariant morphism
$\kappa_0 : U_0 \rightarrow X$, with $X$ endowed with the trivial
$\hat \mu^r$-action,
may be 
embedded equivariantly as an open dense
subset in a 
smooth variety $Z_0$ with good $\hat \mu^r$-action,
endowed with a proper equivariant morphism
$Z_0 \rightarrow X$ extending $\kappa_0$, 
which follows from the appendix of \cite{motivic}
and also from 
Sumihiro's
equivariant completion result \cite{sumi}.
Hence we can 
proceed exactly like in the proof of Theorem \ref{ext} in an
equivariant way, getting existence and unicity of
a $K_0 (\Var_k)$-linear morphism
\begin{equation}\label{equuuuuuubis}
\cS_g : K'_0 (\Var_{X \times \GM^r}^{\GM^r}) \longrightarrow
\cM_{X_0 (g) \times \GM^r \times \GM}^{\GM^r  \times \GM}
\end{equation}
such that, for every 
smooth variety $Z$ with good $\GM^r$-action,
endowed with an equivariant morphism
$p : Z \rightarrow X$
and a monomial morphism ${\mathbf f} : Z \rightarrow \GM^r$
such that
$(p, {\mathbf f}) : Z \rightarrow X \times \GM^r$ is proper
and for every  open
dense subset $U$ of $Z$ which is
stable under the $\GM^r$-action
\begin{equation}\label{skkkkkbis}
\cS_g ([U \rightarrow X \times \GM^r]) = 
\cS_{g \circ p, U}
\end{equation}
in $\cM_{X_0 (g) \times \GM^r \times \GM}^{\GM^r  \times \GM}$.

Let us now prove the compatibility
of the morphism $\cS_g$
with the additional relation (\ref{addrel}).
Let $U$ be a smooth variety over $k$ endowed with a good
$\GM^r$-action, with an equivariant  morphism
$\kappa : U \rightarrow X$,  and with a monomial
morphism $\mathbf{f}_U : U \rightarrow \GM^r$.
Let $q : B \rightarrow U$
be an affine bundle everywhere of rank $s$ with a good affine 
$\GM^r$-action over the action on $U$.
We claim that $U$ may be embedded equivariantly as an open dense
subset in a 
smooth variety $Z$ with good $\GM^r$-action,
endowed with an equivariant morphism
$p : Z \rightarrow X$ extending $\kappa$
and with a monomial morphism ${\mathbf f} : Z \rightarrow \GM^r$
extending ${\mathbf f}_U$,
such that
$(p, {\mathbf f}) : Z \rightarrow X \times \GM^r$ is proper
and  such that, furthermore,
the affine bundle $B$ with its affine $\GM^r$-action
extends
to an affine bundle 
$\widetilde B \rightarrow Z$ with an affine $\GM^r$-action over the action on $Z$
extending the previous one. Indeed, this follows, 
using again the equivalence of categories
of Proposition \ref{equiv} and \ref{twist}, from Lemma \ref{extfib}.
To prove that 
$\cS_g ([B \rightarrow X \times \GM^r])$
does not depend on the affine 
$\GM^r$-action on $B$ over the action on $U$,
it is enough to check that
\begin{equation}\label{yre}
\cS_g ([B \rightarrow X \times \GM^r]) = \LL^s \, \cS_{g} ([U \rightarrow X \times \GM^r]).
\end{equation}
We may assume $(g \circ p)^{-1} (0)$ is 
nowhere dense in $Z$.
Let $h : Y \rightarrow Z$ be a $\GM^r$-log-resolution of
$(Z, (Z \setminus U) \cup (g \circ p)^{-1} (0)$.
We denote by $E_i$, $i$ in $A$, the 
irreducible components of
$h^{-1} ((Z \setminus U) \cup (g \circ p)^{-1} (0))$
and it follows from (\ref{ftyeq}) that,
with the notations of \ref{resolutions} and \ref{eqset},
\begin{equation}\label{ftyeqbis}
\cS_{g} ([U \rightarrow X \times \GM^r]) = - %
\sum_{ \emptyset \not= I \subset C}
(-1)^{\vert I \vert}[U_{I}]
\end{equation}
in $\cM_{X_0
(g)\times \GM^r \times \GM}^{\GM^r \times \GM}$.
Let us consider
the projective bundle $\lambda : Z'  \rightarrow Z$ 
on $Z$, which is the relative projective completion of the bundle
$\widetilde B$.
In particular $Z'$ is endowed with a (projective) $\GM^r$-action.
We consider the pull-back $Y' \rightarrow Y$ of the bundle 
$Z'$ along  the morphism $h$.
We get a proper morphism $h' : Y' \rightarrow Z'$
which is an equivariant log-resolution of
$(Z', (Z' \setminus B) \cup (g \circ p \circ \lambda)^{-1}Ê(0))$.
The set of irreducible components of 
$h'{}^{-1} ((Z' \setminus B) \cup (g \circ p \circ \lambda)^{-1}Ê(0))$
consists of the restriction $E'_i$ of $Y'$ to $E_i$, for $i$ in $A$,
together with $H_{\infty}$, the divisor at infinity of the projective bundle $Y'$.
We set $A' := A \sqcup \{0\}$ and $E'_0 := H_{\infty}$ in order 
to use the notations of \ref{resolutions} and \ref{eqset}
in this setting, adding everywhere  $'$ as an exponent.
In particular
for every non empty subset $I$ of $C'$, we denote by 
$U'_{I}$
the corresponding variety with
$\GM^r\times \GM$-action
and with a monomial morphism 
$(\mathbf{f}'_{I}, g'_{I}) : U'_{I}\rightarrow \GM^r\times \GM$.
Since
$g \circ p \circ \lambda$ is not identically zero on $H_{\infty}$, we have
$C' = C$.
It  follows again from (\ref{ftyeq}) that
\begin{equation}\label{ftyeqter}
\cS_{g} ([B \rightarrow X \times \GM^r]) = - %
\sum_{\emptyset \not= I \subset C}
(-1)^{\vert I \vert}[U'_{I}]
\end{equation}
in $\cM_{X_0
(g)\times \GM^r \times \GM}^{\GM^r \times \GM}$.
Now remark that the natural morphism 
$p_{I} : U'_{I} \rightarrow U_{I}$ 
is an affine bundle of rank $s$,
with an affine $\GM^r \times \GM$-action over the one on 
$U_{I}$.
Furthermore, the 
monomial morphism
$U'_{I}\rightarrow \GM^r\times \GM$
is the composition of 
the 
monomial morphism  
$U_{I}\rightarrow \GM^r\times \GM$ with $p_{I}$.
One deduces that $[U'_{
I}]Ê= \LL^s \, [U_{
I}]Ê$
in $\cM_{X_0 (g) \times \GM^r \times \GM}^{\GM^r  \times \GM}$
and (\ref{yre}) follows.

One then extends
$\cS_g$ by $\cM_k$-linearity to a
$\cM_k$-linear group morphism
\begin{equation}\label{equuuuuuu???}
\cS_g : \cM_{X \times \GM^r}^{\GM^r} \longrightarrow
\cM_{X_0 (g) \times \GM^r \times \GM}^{\GM^r  \times \GM}
\end{equation}
similarly as in the non equivariant case.
\end{proof}

\begin{remark}\label{compbitt}
It follows from our constructions that 
the  morphism
$\cS_g^{\hat \mu^r} : \cM_{X}^{\hat{\mu}^r} \rightarrow
\cM_{X_0(g)}^{\hat{\mu}^{r+1}}$
deduced 
from (\ref{equuuuuuu})
via the canonical
isomorphism (\ref{comp})
is compatible with the one constructed by Bittner
in \cite{bittner2}, modulo the fact  that our additional relation is finer
than hers.
Indeed, they are easily checked to coincide on classes of $\hat \mu^r$-equivariant
morphisms 
$Z \rightarrow X$ with $Z$ smooth and proper.
Note also that diagram \ref{no} is indeed commutative, by construction.
\end{remark}

\begin{lem}\label{extfib}Let $\mathbf{n}$ be in $\NN^r_{>0}$.
Let $X$ be a $k$-variety with trivial $\mu_{\mathbf{n}}$-action
and let
$U$ be a 
smooth
variety with a good $\mu_{\mathbf{n}}$-action
and an equivariant  morphism $\kappa : U \rightarrow X$.
Consider an affine bundle $B \rightarrow U$ with a good
affine $\mu_{\mathbf{n}}$-action over the action on $U$.
Then there exists an equivariant embedding
of $U$ as a dense open set in a smooth
variety $Z$ with good $\mu_{\mathbf{n}}$-action such that
$\kappa$ extends to a proper equivariant morphism
$p : Z \rightarrow X$ and the affine bundle $B$ with its affine $\mu_{\mathbf{n}}$-action
extends
to an affine bundle 
$\widetilde B$ on $Z$ with an affine $\mu_{\mathbf{n}}$-action over the action on $Z$
extending the previous one.
\end{lem}

\begin{proof}Set $G = \mu_{\mathbf{n}}$
and embed $U$ equivariantly in $V$ with a good $G$-action with
$V \rightarrow X$ proper equivariant  extending $\kappa$.
The affine bundle $B  \rightarrow U$ corresponds
to an exact sequence of vector bundles
\begin{equation}\label{bex}
0 \longrightarrow E \longrightarrow F \longrightarrow \cO_U \longrightarrow 0
\end{equation}
on $U$, such that the sheaf of local sections of the affine bundle
is the preimage of $1$ in $F$.
The action of $G$ on $U$ gives a $G$-action
on the exact sequence (\ref{bex}). (By a $G$-action on an $\cO_U$-module
$F$, we mean an isomorphism
$a^* F \rightarrow p^* F$ satisfying the cocycle condition,
with $a : G \times U \rightarrow U$ the action and $p : G \times U \rightarrow U$
the projection on the second factor.)
By blowing up the coherent ideal definining
$V \setminus U$ with the reduced structure we reduce to the case where the inclusion
$j : U \rightarrow V$ is affine. By applying $j_*$ to the
exact sequence (\ref{bex}) and pulling back along $\cO_V \rightarrow j_* \cO_U$,
we extend (\ref{bex}) to an exact sequence of quasi-coherent sheaves with $G$-action
on $V$:
\begin{equation}\label{bextwo}
0 \longrightarrow E' \longrightarrow F'  \longrightarrow \cO_V \longrightarrow 0.
\end{equation}
Let us note that $F'$ is the direct limit of its $G$-invariant coherent subsheaves.
Indeed, this follows from Proposition 15.4 of \cite{lmb}, since (quasi-)coherent sheaves
on the quotient stack $[V / G]$ correspond to 
(quasi-)coherent sheaves with $G$-action on $V$.
It follows 
we may assume the sheaves in (\ref{bextwo}) are coherent.
By restricting to a $G$-stable union of connected components of
$U$,
we may also assume the vector bundle  $E$ is of constant rank $s$ on 
$U$.
Let $q : Z \rightarrow V$ be obtained by taking an equivariant resolution of the 
blow up of the $s$-th Fitting ideal $F_s$ of $E'$,
which is also  the $s+ 1$-th Fitting ideal $F_{s+ 1}$ of $F'$.
Applying $q^*$ to (\ref{bextwo}) and modding out by torsion,
we get an exact sequence of coherent sheaves with $G$-action
\begin{equation}\label{bexter}
0 \longrightarrow \widetilde E \longrightarrow \widetilde F \longrightarrow \cO_{Z} \longrightarrow 0
\end{equation}
on $Z$. 
Let us note that $\widetilde E$ and $\widetilde F$ are in fact locally free. Indeed, 
$Z$ being normal, $\widetilde E$ and $\widetilde F$ are locally free outside a closed subvariety of codimension at least 2, but,  by construction, 
the
Fitting ideals $F_s (\widetilde E)$ and $F_{s + 1} (\widetilde F)$
are invertible,  hence they should be equal to $\cO_Z$.
The preimage of $1$ in $\widetilde F$ is
the sheaf of local sections of an affine bundle with $G$-action
$\widetilde B$ on $Z$ satisfying the required properties.
\end{proof}

\begin{remark}The above proof of Lemma \ref{extfib} 
was explained to us
by Ofer Gabber and works in fact for any linear algebraic group $G$ over $k$.
See also Lemma 7.4 of \cite{bittner} for
a similar, but different, extension lemma.
\end{remark}

\subsection{Compatibility with Hodge realization}
We suppose here that $k = \CC$.
If $X$ is a complex algebraic variety, we denote by
${\rm MHM}_X$ the category of mixed Hodge modules
on $X$, as defined in  \cite{Sa2}. 
We denote by $K_0 ({\rm MHM}_X)$ the corresponding Grothendieck ring.
By addivity, there is a unique 
$\cM_k$-linear
morphism
\begin{equation}
H : \cM_X \longrightarrow K_0 ({\rm MHM}_X)
\end{equation}
such that, for any $p : Z \rightarrow X$  with $Z$ smooth,
$H ([Z])$ is the class of 
the full direct image with compact
supports $Rp_! ({\QQ}_Z)$
in $K_0 ({\rm MHM}_X)$, with ${\QQ}_Z$ the trivial Hodge module on $Z$.
Here we consider $K_0 ({\rm MHM}_X)$
as a  $\cM_k$-module via its $K_0 ({\rm MHM}_{\Spec \CC})$-module structure
and the Hodge realization map $H : \cM_k \rightarrow K_0 ({\rm MHM}_{\Spec \CC})$.
Note that $H (\LL) = [{\QQ}_X (-1)]$.
If $\mu_{{\bf n}} = \mu_{n_1} \times \cdots \times \mu_{n_r}$
acts on $Z$, we may consider
the automorphisms
$T_1$, \dots, $T_r$ on 
the cohomology objects
$R^i p_! ({\QQ}_Z)$
associated respectively to the action of the element with $j$-component
$\exp (2 \pi i / n_j)$ and other components $1$.
If we denote by ${\rm MHM}_X^{r-{\rm mon}}$ the category of mixed Hodge modules on $X$
with $r$ commuting automorphism of finite order,
we get in this way a 
morphism
\begin{equation}\label{tro}
H : \cM_X^{\hat{\mu}^r} \longrightarrow K_0 ({\rm MHM}^{r-{\rm mon}}_X).
\end{equation}
(That the  morphism $H$
is compatible with the additional relation (\ref{addrelmu}), follows from the fact that for
every affine  bundle $p : A \rightarrow Y$ of rank $s$ 
with an affine $\mu_{{\bf n}}$-action
above a $\mu_{{\bf n}}$-action action on $Y$, there is a canonical
equivariant isomorphism $Rp_! ({\QQ}_A) [2s] (s) \simeq \QQ_Y$.)
If $g : X \rightarrow \AA^1$ is a function, there is a nearby cycle functor
$\psi_g : 
{\rm MHM}_X \rightarrow {\rm MHM}_{X_0 (g)}^{\rm mon}$,
cf. \cite{Sa2} \cite{Sa3}, which induces
a morphism
$\psi_g : 
K_0 ({\rm MHM}_X) \rightarrow K_0 ({\rm MHM}_{X_0 (g)}^{\rm mon})$.
By functoriality the  construction extends to 
morphisms
$\psi_g : 
K_0 ({\rm MHM}_X^{r-{\rm mon}}) \rightarrow K_0 ({\rm MHM}_{X_0 (g)}^{r + 1-{\rm mon}})$.

\begin{prop}\label{com}For every $r \geq 0$,
with the notations from Remark \ref{compbitt},
the diagram
\begin{equation}\label{tr}\xymatrix{
\cM_{X}^{\hat{\mu}^r}
\ar[d]^H \ar[r]^{\cS_g^{\hat{\mu}^r}} &
\cM_{X_0 (g)}^{\hat{\mu}^{r + 1}} \ar[d]^H\\
K_0 ({\rm MHM}_X^{r-{\rm mon}})
\ar[r]^{\psi_g}&K_0 ({\rm MHM}_{X_0 (g)}^{r + 1-{\rm mon}})
}
\end{equation}
is commutative.
\end{prop}

\begin{proof}It is enough to prove  that
$H (\cS_g^{\hat{\mu}^r} ([Z \rightarrow X])) = \psi_g (H ([Z \rightarrow X]))$
for $p : Z \rightarrow X$ proper and $Z$ smooth
with   ${\hat{\mu}^r}$-action.  We can further reduce to the case
$(g \circ p)^{- 1}(0)$ is a divisor with normal crossings stable by the $\hat{\mu}^r$-action.
In that case, when $r = 0$,
the statement is proved in \cite{motivic} Theorem 4.2.1 and
Proposition 4.2.3, in a somewhat different language, when $X$ is a point,
but the proof carries over
with no change to general $X$. Since the constructions in loc. cit. may be performed
in an equivariant way in the case of a ${\hat{\mu}^r}$-action, the proof extends directly to the case  $r > 0$.
\end{proof}

\section{Iterated  vanishing cycles}\label{ivcsec}

\subsection{}\label{noname}Let $X$ be a variety endowed 
with the trivial
$\GM^r$-action
and with a function
$g : X \rightarrow \AA^1_k$. Let $U$ be a smooth
$k$-variety of pure dimension $d$ with good $\GM^r$-action
endowed with an equivariant morphism $\kappa : U \rightarrow X$
and with a monomial
morphism $\mathbf{f} = (f_1, \dots, f_r) : U \rightarrow \GM^r$.
Let $U \rightarrow Y$ be an equivariant embedding as a dense open subset of 
a smooth variety $Y$ with a good $\GM^r$-action and with a proper equivariant morphism
$p : Y \rightarrow X$. We assume $(g \circ \kappa)^{-1}Ê(0)$ is nowhere dense in $U$.
Let $h : W \rightarrow Y$ be a 
$\GM^r$-equivariant log-resolution
of $(Y, (Y \setminus U) \cup (g \circ p)^{-1}Ê(0))$.
We shall now explain 
how to compute $\cS_g ([U \rightarrow X \times \GM^r])$
in terms of $W$.
Note that the present set-up is different from the one in Theorem \ref{extequi}.

We denote
by $E_i$, $i$ in $A$, the
irreducible components of
$ h^{-1}((Y \setminus U) \cup (g \circ p)^{-1}Ê(0))$.
We shall use again the notation 
\ref{resolutions} and \ref{eqset}, whenever possible.
Let us assume $I \cap C \not= \emptyset$.
We can still consider the spaces
$U_I$ and the corresponding monomial morphism
$g_I : U_I \rightarrow \GM$.
We denote by $h' : U' \rightarrow U$ the preimage of $U$ in $W$
and we set $F := W  \setminus U'$.
The morphism
$\mathbf{f}: U \rightarrow \GM^r$ extends to a rational map
$\widetilde {\mathbf{f}} : Y  \dashrightarrow (\PP^1_k)^r$.
Furthermore, for $i$ in $A$, there exists integers  $N_i (f_j)$ in $\ZZ$, such that
locally on $W$, each component $\widetilde f_j \circ h$ of $\widetilde {\mathbf{f}} \circ h$ may be written
as
 $u \prod_{i \in A} x_i^{N_i (f_j)}$,
with $u$ a unit, $x_i$ a local equation of
$E_i$. 
Similarly as what we did for $g_I$, 
for every $j$, $1 \leq j \leq r$, we may  define
a rational map
$f_{j, I} : \nu_{E_I} \dashrightarrow \PP^1_k$,
replacing
$N_i (g)$ by $N_i (f_j)$,
and we still denote by $f_{j, I}$ the induced morphism
from $U_I$ to $\GM$. Finally we get
a morphism $\mathbf{f}_I : U_I \rightarrow \GM^r$
which is monomial  for the $\GM^r$-action
by Lemma \ref{promo}.
Similarly as we already observed in \ref{eqset}, this is enough to get that the
morphism
$(\mathbf{f}_{
I}, g_{
I}) : U_{
I}\rightarrow \GM^r\times \GM$
is monomial for the 
$\GM^r\times \GM$-action.
We then denote by  $[U_{
I}]$ the corresponding class in
$\cM_{X_0 (g)\times \GM^r \times \GM}^{\GM^r \times \GM}$.

\begin{lem}\label{promo}Let $W$ be a smooth variety
with a good $\GM^r$-action and let $U$ be a dense open 
stable by the $\GM^r$-action. We assume $F := W \setminus U$ is
a divisor with normal crossings and we denote by $E_i$, $i$ in $A$, its irreducible
components. We consider
a monomial morphism $\mathbf{f} = (f_1, \dots, f_r) : U \rightarrow \GM^r$
and we denote by 
$\widetilde {\mathbf{f}} : W  \dashrightarrow (\PP^1_k)^r$ its extension as a rational
map.
 For any non empty subset $I$ of $A$,
the morphism 
$\mathbf{f}_{
I} : U_{
I}\rightarrow \GM^r$ defined similarly 
as above
is monomial
for the $\GM^r$-action on $U_{
I}$.
\end{lem}

\begin{proof}
Consider the deformation $CW_I$ to the normal cone to $E_I$ in $W$
described in (\ref{desdef}). 
Hence, $CW_{I} := \Spec \cA_{I}$ 
with 
\begin{equation}\label{bgrd}
\cA_{I}:= \sum_{{\bf n} \in \NN^I} \cO_{W \times \AA^I_k}
\left(- \sum_{i \in I} n_i (E_i \times \AA^I_k) \right )
\prod_{i \in I}u_i^{-n_i}.
\end{equation}
Letting $\GM^r$ act trivially on each $u_i$, the $\GM^r$-action on
 $\cO_W$ induces a 
$\GM^r$-action  on $\cA_I$ and on $CW_{I}$.
For $i$ in $A$, we denote by 
$\cJ_i$ the ideal of $\cA_I$ generated by $u_i^{- 1}Ê\cO_W (-E_i)$, resp. $\cO_W (-E_i)$, if
$i \in I$, resp. $i \notin I$, and we set
$\cJ := \prod_{i \in A} \cJ_i$.
We denote by $CW_I^{\circ}$ the complement
in $CW_I$ of the closed subset defined by $\cJ$.
The fiber 
 $CW_I^{\circ} \cap p^{- 1}Ê(0)$ may be identified equivariantly with  $U_{I}$ and
 $CW_I^{\circ}$ with $ p^{- 1}Ê(\GM^I) \simeq U \times \GM^I$ (letting $\GM^r$
 act trivially on $\GM^I$).
 
 On $U \times \GM^I$ we may consider the function
 $(x, u_i) \mapsto f_j (x) \prod_{i \in I} u_i^{- N_i (f_j)}$.
 Similarly as in \ref{desdef}, 
 it extends to a morphism $F_j : CW_I^{\circ} \rightarrow \GM$ whose
 restriction to $U_I$ coincides with $f_{j, I}$.
Let us  consider the  morphism
 $\mathbf{F}  = (F_1, \dots, F_r) :  CW_I^{\circ} \rightarrow \GM^r$.
 Since $\mathbf{f}$ is monomial and $\GM^r$ acts trivially on $u_i$, the restriction of
 $\mathbf{F} $ to the dense open set $U \times \GM^I$
 is monomial, hence $\mathbf{F}$ is monomial and so is its restriction to 
 $U_{ I}$.
 \end{proof}

\subsection{}Let  $\gamma$ and $n$ be in $\NN_{>0}$. 
We keep the notations from
\ref{noname}. In particular  $F = h^{- 1}Ê(Y \setminus U)$.
Let $\varphi$ be in $\cL_{\gamma
n} (W)$ with
$\ord_t
\varphi^*(\cI_F) \leq \gamma n$
 and
$\ord_t g (\varphi) = n$.
 Let $D$ denote the set consisting of all $i$ in $A$ such that
$\varphi (0)$ lies in $E_i$
and consider a local equation $x_i=0$ of
$E_i$ at $\varphi (0)$. By hypothesis, $x_i (\varphi)$ is non zero
in $\cL_{\gamma
n} (\AA^1_k)$, so 
it has a well defined order $\ord_t (x_i (\varphi))$ and angular component
$\ac(x_i (\varphi))$. Writing the component
$\widetilde f_j \circ h$ of $\widetilde {\mathbf{f}} \circ h$
as $u \prod_{i \in D} x_i^{N_i (f_j)}$, with $u$ a unit at $\varphi (0)$,
we set 
\begin{equation}
\ord_t (\widetilde f_j \circ h) (\varphi) := \sum_{i \in D}N_i (f_j) \, \ord_t (x_i (\varphi))
\end{equation} and 
\begin{equation}
\ac (\widetilde f_j \circ h) (\varphi) := u (\varphi (0))  \prod_{i \in D}\ac (x_i (\varphi))^{N_i (f_j)}.
\end{equation}
By abuse of notation, we write $(\widetilde {\mathbf{f}} \circ h)\varphi (0) \in \GM^r$
to mean 
$\ord_t (\widetilde f_j \circ h) (\varphi) = 0$ for every $1 \leq j \leq r$.

Now
we consider the constructible set
\begin{equation}\label{grmter}
\cW_n^{\gamma n}  := \Bigl \{\varphi \in \cL_{\gamma
n} (W)  \Bigm  \vert \ord_t g(\varphi) = n, \ord_t
\varphi^*(\cI_F) \leq \gamma n, (\widetilde {\mathbf{f}} \circ h)(\varphi (0)) \in \GM^r \Bigr\}.
\end{equation}
Similarly as the set in (\ref{grmbis}), $\cW_n^{\gamma n}$ is endowed
with 
a $\GM^r \times \GM$-action and
furthermore the morphism 
$(\ac (\widetilde f_j \circ h), \ac(g)) :
\cW_n^{\gamma n}  \rightarrow \GM^r \times \GM
$
is piecewise monomial.
We denote by 
$[\cW^{\gamma n}_n]$ the corresponding class
in $\cM_{X_0(g)\times \GM^r\times \GM}^{\GM^r\times \GM}$.
Let us consider  the series
\begin{equation}W ^{\gamma} (T) := \sum_{n\geq 1} \, [\cW^{\gamma n}_n]\,\mathbf{L}^{- \gamma  n d}\,T^n
\end{equation} in $\cM_{X_0 (g)\times \GM^r \times \GM}^{\GM^r\times \GM}
[[T]]$.

For $I$ a non empty subset of $A$, 
we consider the cone
\begin{equation}\label{gammai}
\Gamma (I) := \Bigl\{\mathbf{x} \in \RR^I_{>0} \Bigm\vert
\forall j \in \{1, \dots, r\}, \sum_{i \in I} x_i N_i (f_j) = 0\Bigr\}
\end{equation}
and we denote by $d (I)$ its dimension.
We 
shall also consider
the cone
\begin{equation}\label{mgamma}
M_{\gamma} := \Bigl\{\mathbf{x} \in \RR_{>0}^I \Bigm\vert \sum_{i \in I} x_i N_i (\cI_F) \leq \gamma \sum_{i \in I \cap C} x_i N_i (g) \Bigr\}.
\end{equation}

We denote by $\Delta$ the set of non empty subsets $I$ of $A$
such that $\Gamma (I)$ is non empty and is contained in $M_{\gamma}$
for $\gamma \gg 0$.

\begin{prop}\label{formopvar}Let $X$ be a variety with trivial $\GM^r$-action and 
with a function
$g : X \rightarrow \AA^1_k$. 
Let $U$ be a smooth
$k$-variety of pure dimension $d$ with  good $\GM^r$-action
endowed with an equivariant  morphism $\kappa : U \rightarrow X$,
and with a monomial
morphism $\mathbf{f} = (f_1, \dots, f_r) : U \rightarrow \GM^r$.
Let $U \rightarrow Y$ be an equivariant embedding as a dense open subvariety of 
a smooth variety $Y$ with good $\GM^r$-action and with a proper equivariant morphism
$p : Y \rightarrow X$. We assume $(g \circ \kappa)^{-1}Ê(0)$ is nowhere dense in $U$.
Let $h : W \rightarrow Y$ be a 
$\GM^r$-equivariant log-resolution
of $(Y, (Y \setminus U) \cup (g \circ p)^{-1}Ê(0))$.
There exists $\gamma_0$ such that for every $\gamma > \gamma_0$
the series $W^{\gamma}  (T)$ lies
in
$\cM_{X_0 (g) \times \GM}^{\GM} [[T]]_{\rm sr}$
and
$\lim_{T \mapsto \infty} W^{\gamma} (T)$
is independent of
$\gamma > \gamma_0$.
Furthermore, if one sets
$\cW = -\lim_{T \mapsto \infty} W^{\gamma} (T)$, the following holds
\begin{equation}\label{ftyvar}
\cW =  - \sum_{
I \in \Delta
}
(-1)^{d (I)}\,[U_{
I}] 
\end{equation}
in $\cM_{X_0 (g)\times \GM^r \times \GM}^{\GM^r\times \GM}$.
\end{prop}

\begin{proof}Similarly as in the proof of Proposition \ref{formop},
we have
\begin{equation}\label{trtrn}W^{\gamma} (T)  = \sum_{
I \cap C \not= \emptyset
}
[U_{
I}] \, S_{
I} (T)
\end{equation}
with
\begin{equation}
S_{
I} (T) = \sum_{\mathbf{k} \in \Gamma (I) \cap M_{\gamma} \cap \NN^{I}_{ >0} }
\prod_{i \in I} (T^{N_i (g)} \LL^{- 
1 
})^{k_i}.
\end{equation}
The proof now goes on as the proof of Proposition \ref{formop}, 
with $\NN^{I}_{ >0}$ replaced by $\Gamma (I) \cap \NN^{I}_{ >0}$.
Indeed, note that the linear form $\sum_{i \in I \cap C} k_i N_i (g)$ is positive on
$\overline {M_{\gamma}} \setminus \{0\}$
and that $M_\gamma$ is empty if $I\cap
C=\emptyset$. Assume first $I$ lies in $\Delta$ and $I\cap
C\neq \emptyset$. Then
it follows from \ref{chi} that
$\lim_{T \mapsto \infty} S_{
I} (T)  = (-1)^{d (I)}$ for $\gamma \gg 0$. Assume now $I \cap C \not= \emptyset$
and $I \notin \Delta$. In this case, 
necessarily, 
for $\gamma > 0$,
the hyperplane
$\sum_{i \in I} k_i N_i (\cI_F) = \gamma \sum_{i \in I \cap C} k_i N_i (g)$
has a non empty intersection with $\Gamma (I)$. It follows that the Euler characteristic
$\chi (\Gamma (I) \cap M_{\gamma})$ is equal to zero.
\end{proof}

\begin{prop}\label{eqstrict}Let $X$ be a variety with trivial $\GM^r$-action and 
with a function
$g : X \rightarrow \AA^1_k$. 
Let $U$ be a smooth
$k$-variety of pure dimension $d$ with  good $\GM^r$-action
endowed with an equivariant  morphism $\kappa : U \rightarrow X$,
and with a monomial
morphism $\mathbf{f} = (f_1, \dots, f_r) : U \rightarrow \GM^r$.
Let $U \rightarrow Y$ be an equivariant embedding as a dense open subvariety of 
a smooth variety $Y$ with good $\GM^r$-action and with a proper equivariant morphism
$p : Y \rightarrow X$. We assume $(g \circ \kappa)^{-1}Ê(0)$ is nowhere dense in $U$.
Let $h : W \rightarrow Y$ be a 
$\GM^r$-equivariant log-resolution
of $(Y, (Y \setminus U) \cup (g \circ p)^{-1}Ê(0))$.
Then, with the previous notation, we
have
\begin{equation}\label{pint}
\cS_{g} ([U\rightarrow X \times \GM^r])= - \sum_{
I \in \Delta
}
(-1)^{d (I)}\,[U_{
I}] 
\end{equation}
in $\cM_{X_0 (g)\times \GM^r \times \GM}^{\GM^r\times \GM}$.
\end{prop}

\begin{proof}We may reduce to the case where
the morphism
$\mathbf{f}: U \rightarrow \GM^r$ extends to a morphism
$\widetilde {\mathbf{f}} : Y  \rightarrow (\PP^1_k)^r$.
Indeed, there exists an equivariant embedding
$U \rightarrow Y'$  of $U$ as a dense open subvariety of 
a smooth variety $Y'$ with a good $\GM^r$-action and with a proper equivariant morphism
$p' : Y' \rightarrow X$ such that
$\mathbf{f}$ extends to a morphism
$\widetilde {\mathbf{f}}': Y' \rightarrow (\PP^1_k)^r$.
We may furthermore assume there is a $\GM^r$-equivariant proper morphism
$Y' \rightarrow Y$ which restricts to the identity on $U$.
Let $h' : W' \rightarrow Y'$ be a 
$\GM^r$-equivariant log-resolution
of $(Y', (Y' \setminus U) \cup (g \circ p')^{-1}Ê(0))$.
We may also assume there is a $\GM^r$-equivariant proper morphism
$W' \rightarrow W$ such that the diagram
\begin{equation}\xymatrix{
W'
\ar[d] \ar[r]^{h'} &
Y'
 \ar[d]\\
W
\ar[r]^{h}&Y
}
\end{equation}
is commutative.

Consider $\cW'$ defined as $\cW$ but using $W'$ instead of $W$.
Since, temporarily, we shall work on $W'$ and not on $W$,
we denote by $E_i$, $i$ in $A$, the
irreducible components of
$ h'{}^{-1}((Y' \setminus U) \cup (g \circ p')^{-1}Ê(0))$, and keep the
previous notation, but for 
$W'$ instead of $W$.
We have 
\begin{equation}
\cW' = - \lim_{T \mapsto \infty}Ê\sum_{I \cap C \not=0}Ê[U_I] S_I (T),
\end{equation}
while, computing $W^{\gamma}Ê(T)$ on $W'$ using the change of
variable formula, or more precisely Lemma 3.4 in
\cite{arcs}, one gets
\begin{equation}
\cW = - \lim_{T \mapsto \infty}Ê\sum_{I \cap C \not=0}Ê[U_I] \widetilde S_I (T)
\end{equation}
with 
\begin{equation}
\widetilde S_{I} (T) = \sum_{\mathbf{k} \in \Gamma (I) \cap M_{\gamma} \cap \NN^{I}_{ >0} }
\prod_{i \in I} (T^{N_i (g)} \LL^{- m_i })^{k_i},
\end{equation}
with $m_i \geq 1$. It follows that  $\cW' = \cW$,
and by Proposition \ref{formopvar} we can assume $Y = Y'$ and $W = W'$.

Consider $Z := (\widetilde {\mathbf{f}} \circ h)^{-1}(\GM^r)$ in $W$.
Note that the image of the morphism
$Z \rightarrow W \times \GM^r$ given by the inclusion on the first factor
and by the restriction of $\widetilde {\mathbf{f}} \circ h$ on the second factor
is the closure of the image of the inclusion $U' \rightarrow W \times \GM^r$,
with $U'$ the preimage of $U$ in $W$.
It follows that the morphism $(q, \widetilde {\mathbf{f}} \circ h_{\vert Z}):
Z \rightarrow X \times \GM^r$ given by composition
with $p \circ h$ on the first factor is proper.
Since $Z$  is smooth and the morphism
$Z \rightarrow X \times \GM^r$ is proper, it follows from 
(\ref{skkkkk}) that 
\begin{equation}\label{hpint}
\cS_g ([U  \rightarrow X \times \GM^r]) = 
\cS_{g \circ q, U'}
\end{equation}
in $\cM_{X_0 (g) \times \GM^r \times \GM}^{\GM^r  \times \GM}$.
Note also that, since $\mathbf{f}: U \rightarrow \GM^r$ extends to a morphism
$\widetilde {\mathbf{f}} : Y  \rightarrow (\PP^1_k)^r$, 
for a  
subset $I$ of $A$ with
$I \cap C \not= \emptyset$,
$\Gamma (I)$ is non empty if and only if
$E_I^{\circ}$ is contained in  $(\widetilde {\mathbf{f}} \circ h)^{-1}(\GM^r)$.
Furthermore if these conditions hold, $\Gamma (I) = \RR^I_{>0}$.
It follows that $\Delta$ consists exactly of those non empty subsets of $C$
for which $E_I^{\circ}$ is contained in  $(\widetilde {\mathbf{f}} \circ h)^{-1}(\GM^r)$,
hence the right hand side of
(\ref{pint}) may be rewritten
as
\begin{equation}
-
\sum_{\sur{
\emptyset \not= I \subset C
}{E_I^{\circ} \subset (\widetilde {\mathbf{f}} \circ h)^{-1}(\GM^r)
}}
(-1)^{\vert I \vert}\,[U_{
I}],
\end{equation}
and 
the required equation (\ref{pint}) follows now from (\ref{hpint}) and 
(\ref{ftyeq}).
\end{proof}

\subsection{Iterated vanishing cycles}\label{resol}
Now we consider a smooth variety $X$ of pure dimension $d$
with  two functions $f : X \rightarrow  \AA^1_k$ and $g : X
\rightarrow  \AA^1_k$. The motivic Milnor fiber $\cS_f$
lies in $\cM_{X_0 (f) \times \GM}^{\GM}$. We still denote
by $g$ the function $X_0 (f) \times \GM \rightarrow
\AA^1_k$ obtained by composition of $g$ with the projection
$X_0 (f) \times \GM \rightarrow X$. Hence, thanks to
\ref{eqset}, we may consider the image
\begin{equation}\cS_g (\cS_f) = \cS_g (\cS_f ([X \rightarrow X]))
\end{equation}
of 
$\cS_f = \cS_f ([X \rightarrow X])$
by the nearby cycles
morphism 
\begin{equation}\cS_g : \cM_{X_0 (f)\times \GM}^{\GM} \longrightarrow 
\cM_{(X_0 (f) \cap X_0
(g))\times \GM^2}^{\GM^2}\end{equation}
 which lies in $\cM_{(X_0 (f) \cap X_0
(g))\times \GM^2}^{\GM^2}$.

We shall now give an explicit
description of $\cS_g (\cS_f)$ in terms of a log-resolution  $h : Y \rightarrow X$
of  $(X, X_0 (f) \cup X_0 (g))$.
We shall denote by 
$E_i$, $i$ in $A$,
the irreducible components of
$h^{-1} (X_0 (f) \cup X_0 (g))$ and we
shall consider the sets
\begin{equation}
B = \Bigl\{i \Bigm \vert N_i (f) > 0\Bigr\} \quad
\text{and} \quad
C = \Bigl\{i \Bigm \vert N_i (g) > 0\Bigr\}.
\end{equation}
Recall, cf. \ref{resolutions},  that we denoted by $U_I^J$,
for $J \subset I$, the fiber product of the restrictions of
the $\GM$-bundles $U_{E_i}$, for $i$ in $J$, to
$E_I^{\circ}$. Assume $J := I \cap C$ and $K := I \setminus
C$ are both non empty. We now consider the fiber product
$U_{K,J} := U_I^K \times_{E_I^{\circ}} U_I^J$ which has the
same underlying variety than $U_I$. There is a natural
$\GM^2$-action on $U_{K,J}$, the first, resp. second,
$\GM$-action being the diagonal action on $U_I^K$, resp.
$U_I^J$, and the trivial one on the other factor.
The  morphism
$(f_I,g_I) : U_I = U_{K,J} \rightarrow  \GM^2$
being monomial,
the morphism
$(h\circ \pi_I, f_I,g_I) : U_{I} \rightarrow (X_0 (f) \cap
X_0 (g)) \times \GM^2$ has a class
in $\cM_{(X_0 (f) \cap X_0 (g)) \times
\GM^2}^{\GM^2}$ that  we denote by $[U_{K,J}]$.

\begin{theorem}\label{reee}With the previous notations, we have
\begin{equation}\label{sgsf}
\cS_g (\cS_f) =\sum_{\sur{I\cap
C=J\neq \emptyset}{I\setminus C=K\neq \emptyset}}(-1)^{|I|}\,[U_{K,J}]
\end{equation}
in $\cM_{(X_0 (f) \cap X_0 (g)) \times \GM^2}^{\GM^2}$.
\end{theorem}

\begin{proof}Consider the inclusions
$i : X_0 (g) \times \GM \hookrightarrow X\times \GM$ and $j
: (X \setminus X_0 (g)) \times \GM \hookrightarrow X\times
\GM$. Note that $\cS_f - j_! (\cS_{(f_{\vert X \setminus
X_0 (g)})})$ is supported by $X_0 (g) \times \GM$, that is,
is of the form $i_! (\cA)$. Since $Y \setminus Y_0 (g \circ
h)$ is a log-resolution of $(X \setminus X_0 (g), X_0 (f)
\setminus X_0 (g))$,
one deduces from the proof of 
(\ref{dl2})  that
\begin{equation}\label{vill}
\cS_{f}  - \Bigl(-\sum_{\sur{K \cap C = \emptyset}{K \not= \emptyset}}\,
(-1)^{|K|}\,[U_K \rightarrow X_0 (f) \times \GM]\Bigr)
\end{equation}
is supported by $X_0 (g) \times \GM$, hence, since
$\cS_g$ is zero on objects of the form $i_! (\cA)$,
we deduce that
\begin{equation}\label{ccc}\cS_g (\cS_{f}) = \cS_g
 \Bigl( -\sum_{\sur{K \cap C = \emptyset}{K \not= \emptyset}}\,  (-1)^{|K|}\,[U_K
\rightarrow X_0 (f) \times \GM]\Bigr).
\end{equation}

To conclude it is enough to check the following equality in
$\cM_{(X_0 (f) \cap X_0 (g)) \times \GM^2}^{\GM^2}$, for every non empty 
subset $K$ of $A$ such that $K\cap C = \emptyset$:
\begin{equation}\label{guio}
\cS_{g} ( [U_K \rightarrow X_0 (f) \times \GM]) = - \sum_{\emptyset \not= J \subset C}
(-1)^{\vert J \vert} [U_{K, J}].
\end{equation}
This will  follow from Proposition \ref{eqstrict}. Indeed,
let us consider the projective bundle
$\pi_K : \overline{\nu}_{E_K} \rightarrow E_K$
with the  $\GM$-action extending the diagonal one on $\nu_{E_K}$.
Let us set $A' := A\sqcup \{\infty\}$.
The complement of $U_K$ in
$ \overline{\nu}_{E_K}$ is a divisor with normal crossings
whose irreducible components are:
\begin{enumerate}
\item[-]the divisors $E'_j := \pi_K^{-1} (E_{K \cup \{j\}})$, for $j \notin K$ such that
$E_{K \cup \{j\}} \not= \emptyset$
\item[-] the divisor at infinity $E'_{\infty} :=\overline{\nu}_{E_K}\setminus \nu_{E_K}$
\item[-] the divisors $E'_{i}$, for $i$ in $K$, defined as the closure of
  the fiber product, above $E_K$, of the zero section of $\nu_{E_i}$
  with the $\nu_{E_\ell}$, $\ell$ in $K$, $\ell\neq i$.
\end{enumerate}
Note that all these divisors are stable by the $\GM$-action.
We shall use the notations of \ref{resolutions} and \ref{eqset} with an exponent $'$.

We now determine the set 
$\Delta$ of non empty subsets $J'$ of $A'$
such that $\Gamma (J')$ is non empty and is contained in $M_{\gamma}$
for $\gamma \gg 0$, with the notation of  (\ref{gammai}) and
(\ref{mgamma}).

Note that for $\Gamma (J')$ to be  non empty it is necessary that if
$N_i(f_K) > 0$ (resp. $N_i(f_K) < 0$)
for some $i$ in $J'$, then for some $i'$ in $J'$, 
$N_{i'}(f_K) < 0$ (resp. $N_{i'}Ê(f_K) > 0$).
This forces $J'$ to be either of the form $J \sqcup \{\infty\}$
with $J \cap B \not= \emptyset$ or of the form $J$ with $J \cap B= \emptyset$.
In each case, the condition that
$\Gamma (J')$ is contained in $M_{\gamma}$
for $\gamma \gg 0$ implies that $J \subset C$ and furthermore that
$d (J') = \vert J \vert$.
We deduce that
$J \sqcup \{\infty\}$
belongs to $\Delta$ if and only $J \cap B \not=\emptyset$ and $J \subset C$  and 
that $J$
belongs to $\Delta$ if and only $J \cap B =\emptyset$, $J \subset C$
and $J \not= \emptyset$.
It follows from
Proposition \ref{eqstrict}
that
\begin{equation}\label{spint}
\cS_{g} ([U_K \rightarrow X_0 (f)
\times \GM])= - \sum_{\sur{\emptyset \not= J \subset C}{J \cap B = \emptyset}} 
(-1)^{\vert J \vert}\,[U'_{J}] -   \sum_{\sur{ J \subset C}{J \cap B  \not= \emptyset}}
(-1)^{\vert J \vert}\,[U'_{J \sqcup \{\infty\}}]
\end{equation}
in $\cM_{X_0 (g)\times \GM^r \times \GM}^{\GM^r\times \GM}$.
To conclude it is enough to note that if
$\emptyset \not= J\subset C$ and $J \cap B = \emptyset$,
then $[U'_{J}]  = [U_{K, J}]$ and that if $ J \subset C$
and $J \cap B  \not= \emptyset $ we have $ [U'_{J \sqcup
\{\infty\}}]  = [U_{K, J}]$.


Let us prove the second equality. We consider the image $\PP(U_{E_K})$
of  $U_{E_K}$ in
$\PP(\nu_{E_K})$ and note that
the canonical morphism
$U_{E_K}\rightarrow \PP(U_{E_K})$ is a $\GM$-bundle, namely the 
restriction to $U_{E_K}$ of the tautological  line bundle  on $ \PP(\nu_{E_K})$.
We
identify $E''_\infty :=
E'_\infty \setminus \cup_{i \in K} E'_i$
with  $\PP(U_{E_K})$. 
The restriction of the tautological line bundle
to $\PP(U_{E_K})$ is dual to the restriction to
$E''_\infty$ of the normal bundle to $E'_\infty$ in
$\overline{\nu}_{E_K}$. We have now two $\GM$-bundles on
$E''_\infty=\PP(U_{E_K})$, namely $U_{E_K}$ and the
restriction, we shall denote by $U''_{\infty}$,  to $E''_\infty$ of the complement
$U_{E'_\infty}$ of the zero section in the normal bundle
$\nu_{E'_\infty}$. 
Let us denote by $U_{E_K}^a$ the $\GM$-bundle
$U_{E_K}$ endowed with the inverse $\GM$-action.
The antipody $a : U_{E_K} \rightarrow U_{E_K}^a$ whose restriction to the fibers is given by
$t \mapsto t^{-1}$ is an isomorphism of $\GM$-bundles with $\GM$-action.
By the above description, $U''_{\infty}$ may be identified, as a 
$\GM$-bundle with $\GM$-action, with
the $\GM$-bundle $U_{E_K}^a$.

Now consider the function $f_K$ on
$\nu_{E_K}$. It induces a rational map $\widetilde f_K$ on
$\overline{\nu}_{E_K}$.
Let us check that under the above isomorphism, the restriction
$f_K : U_{E_K} \rightarrow \AA^1_k$ composed with the automorphism $a$
corresponds to the morphism
$f'_{\infty} : U''_{\infty} \rightarrow \AA^1_k$
obtained from $\widetilde f_K$ by the construction of \ref{noname}.
Indeed, let 
$U$ be an
open subset of $E_{K}$ above which the bundle $\nu_{E_K}$ is
trivial, isomorphic to $U\times \AA_k^{K}$. 
Denote by $w_i$, for $i$ in $K$, the coordinates on $\AA_k^{K}$.
Fix $\ell$ in $K$. 
The restriction of $U_{E_K}$ to $U$ may be identified, equivariantly, with
$U \times \PP (\GM^K) \times \GM$
by
$(u, (w_i)_{i \in K}) \mapsto (u, (x_i  = \frac{w_i}{w_{\ell}})_{i \in K \setminus \{\ell\}}, t = w_{\ell})$, where $\GM$ acts trivially on the first two factors and by multiplicative
translation on the last one, with $(x_i)_{i \in K \setminus \{\ell\}}$ the standard coordinates on
$\PP (\GM^K) \simeq \GM^{K \setminus \{\ell\}}$.
If 
the restriction of $f_K$ to
$U\times \GM^{K} $ is given by
$v(u)\prod_{i\in K} w_i^{N_i}$, it may be rewritten, under the above identification,  as
$v(u)\prod_{i\in K \setminus \{\ell\}} x_i^{N_i} t^{\sum_{i\in K}N_i}$.
Composing with the antipody $a$
we get the function
$v(u)\prod_{i\in K \setminus \{\ell\}} x_i^{N_i} t^{- \sum_{i\in K}N_i}$
which corresponds to  the  restriction of the function $f'_\infty$ on the corresponding
open subset.

If $J$ is a subset of $C$ such that $E_{K\sqcup J}\neq
\emptyset$, 
it follows from the ``transitivity'' property
described in
 \ref{desdef} that $f_{K\sqcup J}$ can be retrieved directly from $f_K:
\nu_{E_K}\rightarrow \AA_k^1$
and similarly, the rational map
$f'_{J\sqcup \{\infty\}}$ can be retrieved directly from
the rational map $f'_\infty$ (obtained from $\widetilde f_K$ by the construction of \ref{noname})
on $\nu_{E'_\infty}$. It follows that, under the isomorphism 
 between $U_{K\sqcup J}$ and  $U'_{J\sqcup \{\infty\}}$ induced by $\varphi$,
$f_{K\sqcup
J}$ corresponds to $f'_{J\sqcup \{\infty\}}$.The same argument works for the functions induced by $g$ on $U'_{K\sqcup J}$
and $U_{J\sqcup \{\infty\}}$ (note that in fact $N_i(g) = 0$  for all $i$ in $K$ and
 $N_\infty(g) = 0$).

The first equality, which is easier,  is  checked similarly using
$E'_J=\pi_K^{-1}(E_{K \cup J})$
and the canonical isomorphism of bundles $\nu_{E'_J} \simeq
(\pi_{K \vert E'_J})^* (\nu_{E_{J}\vert E_{K \cup J}})$, for $J \subset C$.
\end{proof}

\section{Convolution and the main result}\label{cmr}
\subsection{Convolution}\label{cov}
Let us denote by $a$ and $b$ the coordinates on each  factor
of $\GM^2$.
Let $X$ be a variety.
We denote by
$i : X \times (a + b)^{-1} (0) \rightarrow X \times \GM^2$
the inclusion of the antidiagonal and by $j$ the inclusion of its complement.
We consider the morphism
\begin{equation}
 a+b  :   X \times \GM^2\setminus (a+b)^{-1}(0)
 \longrightarrow X \times \GM
\end{equation}
which is the identity on the $X$-factor and is equal to
$a + b$ on the $\GM^2\setminus (a+b)^{-1}(0)$-factor.
We denote by $\pr_1$ and $\pr_2$
the projection of $X \times \GM\times (a+b)^{-1}(0)$ on
$X \times \GM$ and $X \times (a+b)^{-1}(0)$,
respectively.

If  $A$ is an object in $\cM_{X \times \GM^2}$,
the object
\begin{equation}
\Psi^0_{\Sigma} (A) := - (a + b)_! j^* (A) + \pr_{1!} \pr_2^* i^*(A)
\end{equation}
lives in $\cM_{X \times \GM}$.
We now explain how to
lift
$\Psi^0_{\Sigma}$ to a $\cM_k$-linear group morphism
$\Psi_{\Sigma} : \cM_{X \times \GM^2}^{\GM^2} \rightarrow
\cM_{X \times \GM}^{\GM}$.

Let $A$ be  an object in
$\Var_{X \times \GM^2}^{\GM^2, (n, m)}$
with class $[A]$ in
$\cM_{X \times \GM^2}^{\GM^2, (n, m)}$.
It is endowed 
with a $\GM^2$-action $\alpha$
for which the morphism to $\GM^2$ is diagonally monomial of weight $(n, m)$.
We may consider the $\GM$-action $\widetilde \alpha$ on $A$
given by $\widetilde \alpha (\lambda) x = \alpha (\lambda^m, \lambda^n) x$.
With some obvious  abuse of notations,
$(a + b)_! j^* ([A])$
is the class of 
$a + b : A_{\vert a + b \not= 0} \rightarrow X \times \GM$.
If we endow $A_{\vert a + b \not= 0}$ with the $\GM$-action induced by $\widetilde \alpha$,
the morphism $( a + b) : A_{\vert a + b \not= 0} \rightarrow  \GM$ is diagonally
monomial of weight $nm$,
The   term
$\pr_{1!} \pr_2^* i^*([A])$
is the class of 
$ A_{\vert a + b = 0}  \times \GM \rightarrow X \times \GM$,
the morphism to $\GM$ being  the projection
on the $\GM$-factor. We endow 
$ A_{\vert a + b = 0}  \times \GM$ with the $\GM$-action induced by $\widetilde \alpha$
on the first factor and the action
$(\lambda, z) \mapsto \lambda^{nm}z$ on the second factor.
Hence we may set
$\Psi_{\Sigma}^{n, m} ([A]) = - [a + b : A_{\vert a + b \not= 0} \rightarrow X \times \GM]  +
[A_{\vert a + b = 0}  \times \GM \rightarrow X \times \GM]
$
in 
$\cM_{X \times \GM}^{\GM, nm}$ and extend
this construction in  a unique 
way to a 
$\cM_k$-linear group morphism
\begin{equation}
\Psi_{\Sigma}^{n, m} : \cM_{X \times \GM^2}^{\GM^2, (n, m)} \longrightarrow
\cM_{X \times \GM}^{\GM, nm}.
\end{equation}
These morphisms being compatible with the morphisms induced
by the transition morphisms
of (\ref{dlel}),  we get after passing to the colimit a 
$\cM_k$-linear group morphism
\begin{equation}
\Psi_{\Sigma} : \cM_{X \times \GM^2}^{\GM^2} \longrightarrow
\cM_{X \times \GM}^{\GM}.
\end{equation}

Let us now explain the relation of
$\Psi_{\Sigma}$
with the convolution product as considered in
\cite{ts}, \cite{looi} and \cite{barc}.
There is a canonical
morphism
\begin{equation}\label{noto}
\cM_{X \times \GM}^{\GM}
\times
\cM_{X \times \GM}^{\GM}
\longrightarrow
\cM_{X \times \GM^2}^{\GM^2}
\end{equation}
sending $(A, B)$ to
$A \boxtimes B$, the fiber product over $X$ of $A$ and $B$,
therefore we may define
\begin{equation}
\ast :
\cM_{X \times \GM}^{\GM}
\times
\cM_{X \times \GM}^{\GM}
\longrightarrow
\cM_{X \times \GM}^{\GM}
\end{equation}
by
\begin{equation}
A \ast B = \Psi_{\Sigma} (A \boxtimes B).
\end{equation}

If $S$ is in  $\Var_X^{\mu_n}$, resp. in
$\Var_X^{\mu_n \times \mu_n}$,   we  denote by $[S]$ the
corresponding  class in
$\cM_{X \times \GM}^{\GM}$, resp. in $\cM_{X \times \GM^2}^{\GM^2}$,
via the isomorphism
(\ref{comp}).
Consider the
Fermat curves
$F_1^n$ and $F_0^n$ 
defined respectively
by $x^n + y^n = 1$
and $x^n + y^n = 0$ in $\GM^2$ with their  standard
$\mu_n \times \mu_n$-action.
If  $A$ is a variety in 
$\Var_X^{\mu_n \times \mu_n}$, we have 
\begin{equation}\label{prefermat}
\Psi_{\Sigma} 
([A]) = - [F_1^n \times^{\mu_n \times \mu_n} A]
+
[F_0^n \times^{\mu_n \times \mu_n} A],
\end{equation}
the $\mu_n$-action
on each term in the right hand side of
(\ref{prefermat}) being the diagonal one.
In particular, if 
$A$ and $B$ are two varieties in
$\Var_X^{\mu_n}$, the convolution product $[A] \ast [B]$
is given by
\begin{equation}\label{fermat}
[A] \ast [B] = - [F_1^n \times^{\mu_n \times \mu_n} (A \times_X B)]
+
[F_0^n \times^{\mu_n \times \mu_n} (A \times_X  B)].
\end{equation}
The convolution product 
in \cite{looi} and
\cite{barc} was defined when $k$ contains all roots of unity.
Since as soon as $k$ contains a $n$-th root of
$-1$ we have
$[F_0^n \times^{\mu_n \times \mu_n} (A \times_X B)] = (\LL - 1) [(A \times_X B) / \mu_n]$,
one gets that the convolution product
in \cite{looi} and
\cite{barc}, when defined, coincides with the one in (\ref{fermat}).

\begin{prop}\label{ass}The convolution product
on $\cM_{X \times \GM}^{\GM}$ is commutative and associative.
The unit element for the convolution product
is $1$, the class of the identity $X \times \GM \rightarrow
X \times \GM$ with  the standard $\GM$-action on the
$\GM$-factor.
\end{prop}

\begin{proof}Commutativity being clear, let us prove the statement concerning
associativity and unit element.
For simplicity of notation
we shall assume $X$ is a point and
we shall first ignore the $\GM$-actions, that is we shall
prove the corresponding statements for
$\cM_{\GM}$.
Consider $a : A \rightarrow \GM$,
$b : B \rightarrow \GM$, $c : C \rightarrow \GM$.
By definition the convolution product
$A \ast B$ (with some  abuse of notation, we shall denote by the same
symbol varieties over $\GM$ and their class in $\cM_{\GM}$)
is equal to
\begin{equation}
-[ a+ b : (A \times B)_{\vert a + b \not= 0} \rightarrow \GM]
+ [z : (A \times B \times \GM)_{\vert a + b = 0} \rightarrow \GM],
\end{equation}
with $z$ the standard coordinate
on $\GM$.

Associativity follows from the following claim:
$(A \ast B) \ast C$ is equal to \begin{equation}\label{claim}
[ a+ b + c : (A \times B \times C)_{\vert a + b+ c \not= 0} \rightarrow \GM]
- [z : (A \times B \times C \times \GM)_{\vert a + b + c= 0} \rightarrow \GM].
\end{equation}
Indeed, $(A \ast B) \ast C$ may be written as a sum of four terms.
The first one,
\begin{equation}
[ a+ b + c : (A \times B \times C)_{\vert \sur{a + b+ c \not= 0}{a + b \not=0}} \rightarrow \GM]
\end{equation}
may be rewritten as
\begin{equation}\label{eq1}
[ a+ b + c : (A \times B \times C)_{\vert  a + b+ c \not= 0}
\rightarrow \GM]
- [ c : (A \times B \times C)_{\vert a + b =0}
\rightarrow \GM].
\end{equation}
The second one,
\begin{equation}
-
[z : (A \times B \times C \times \GM)_{\vert \sur{a + b + c= 0}{a + b \not=0}}
\rightarrow \GM]
\end{equation}
may be rewritten as
\begin{equation}\label{eq2}
- [z : (A \times B \times C \times \GM)_{\vert a + b + c= 0}
\rightarrow \GM].
\end{equation}
The third one
\begin{equation}
- [c  + z : (A \times B \times C \times \GM)_{\vert \sur{a + b = 0}{c + z \not=0}}
\rightarrow \GM]
\end{equation}
may be rewritten as
\begin{equation}\label{eq3}
- [u : (A \times B \times C \times \GM)_{\vert \sur{a + b = 0}{u \not=c}}
\rightarrow \GM],
\end{equation}
since the corresponding spaces are isomorphic via
$(\alpha, \beta, \gamma, z) \mapsto (\alpha, \beta, \gamma,  u = c (\gamma) + z)$.
Here $u$ is a coordinate on some other copy of $\GM$.
The fourth term,
\begin{equation}
[u : (A \times B \times C \times \GM \times \GM)_{\vert \sur{a + b = 0}{c + z = 0}}
\rightarrow \GM]
\end{equation}
may be rewritten as
\begin{equation}\label{eq4}
[u : (A \times B \times C \times \GM)_{\vert a + b = 0}
\rightarrow \GM].
\end{equation}
One deduces (\ref{claim}) by summing up
(\ref{eq1}), (\ref{eq2}), (\ref{eq3}), and (\ref{eq4}).

For the statement concerning
the unit element, one writes
$A \ast \GM$ as
\begin{equation}
- [a + z : (A \times \GM)_{\vert a + z \not= 0} \rightarrow \GM]
+ [ u : (A \times \GM \times \GM)_{\vert a + z = 0} \rightarrow \GM].
\end{equation}
Since the first term may be rewritten as
\begin{equation}
- [u : (A \times \GM)_{\vert a \not= u} \rightarrow \GM]
\end{equation}
and the second term as
\begin{equation}
[u : (A \times \GM) \rightarrow \GM],
\end{equation}
it follows that $A \ast \GM$ is equal to (the class of) $A$
in $\cM_{\GM}$.
The proofs for  general $X$ are just the same.
As for $\GM$-actions, since  by the very constructions they
are diagonally monomial of the same weight on each factor,
all identifications we made
are compatible with the $\GM$-actions,
and all
statements still hold in
$\cM_{X \times \GM}^{\GM}$.
\end{proof}

\begin{remark}
Proposition \ref{ass}, modulo
the isomorphism (\ref{comp}),
is already stated in  \cite{barc}.
\end{remark}

\subsection{}In fact, associativity already holds
at the $\Psi_{\Sigma}$-level.
To formulate this, we need to introduce some more notation.

Let us denote by $a$, $b$ and $c$ the coordinates on each  factor
of $\GM^3$.
For $X$ a variety,
we denote by
$i$ the inclusion
$X \times (a + b + c)^{-1} (0) \hookrightarrow X \times \GM^3$
and by $j$ the inclusion of the complement.
We consider the morphism
\begin{equation}
 a+b +c  :   X \times \GM^3\setminus (a+b + c)^{-1}(0)
 \longrightarrow X \times \GM
\end{equation}
which is the identity on the $X$-factor and is equal to
$a + b + c$ on the $\GM^3\setminus (a+b + c)^{-1}(0)$-factor.
We denote by $\pr_1$ and $\pr_2$
the projection of $X \times \GM\times (a+b+ c)^{-1}(0)$ on
$X \times \GM$ and $X \times (a+b+ c)^{-1}(0)$,
respectively.

If  $A$ is an object  in $\cM_{X \times \GM^3}^{\GM^3}$,
we consider the object
\begin{equation}
\Psi^0_{\Sigma_{123}} (A) :=  (a + b + c)_! j^* (A) - \pr_{1!} \pr_2^* i^*(A),
\end{equation}
in $\cM_{X \times \GM}$. Similarly as in \ref{cov} we extend
$\Psi^0_{\Sigma_{123}}$ to a $\cM_k$-linear group morphism
$\Psi_{\Sigma_{123}} : \cM_{X \times \GM^3}^{\GM^3} \rightarrow \cM_{X \times \GM}^{\GM}$.
We
 denote by $A_{ij}$ the object $A$ viewed
as an element in  $\cM_{X \times \GM^2}^{\GM^2}$
by forgetting the projection and the action corresponding
to the $k$-th $\GM$-factor, with $\{i, j, k\} = \{1, 2, 3\}$.
The object
$\Psi_{\Sigma} (A_{ij})$
may now be endowed with a second projection to $\GM$ and a second
$\GM$-action, namely those corresponding to
the $k$-th $\GM$-factor, so we get in fact an element in
$\cM_{X \times \GM^2}^{\GM^2}$ we denote by
$\Psi_{\Sigma_{ij}} (A)$.

\begin{prop}\label{superass}Let
$A$ be an object  in $\cM_{X \times \GM^3}^{\GM^3}$. For every
$1 \leq i < j \leq 3$, we have
\begin{equation}
\Psi_{\Sigma_{123}} (A) = \Psi_{\Sigma} ( \Psi_{\Sigma_{ij}} (A)).
\end{equation}
\end{prop}

\begin{proof}The proof is the same as the one for
associativity in
Proposition \ref{ass}. Indeed, one just has to replace everywhere
$A \times B \times C$ by $A$ in the proof,
and to remark that (\ref{claim}) then becomes nothing else than
$\Psi_{\Sigma_{123}} (A)$.
\end{proof}

\subsection{}Let us consider again a smooth variety
$X$ of pure dimension $d$ with two functions $f$ and $g$ from $X$ to $\AA^1_k$.
Let us denote by
$i_1$ and $i_2$ the inclusion of
$(X_0 (f) \cap X_0 (g)) \times \GM$
in $X_0 (f) \times \GM$ and
$X_0 (f + g^N) \times \GM$, respectively.

We can now state the main result of this paper.

\begin{theorem}\label{MT}Let $X$ be a smooth variety of pure dimension $d$,
and $f$ and $g$ be two functions  from $X$ to $\AA^1_k$.
For every $N >  \gamma ((f), (g))$, the equality
\begin{equation}\label{mt}
i_1^* \cS_{f}^\phi - i_2^*\cS_{f+g^N}^\phi = \Psi_{\Sigma}
(\cS_{g^N}(\cS_f^\phi))
\end{equation}
holds in $\cM_{(X_0 (f) \cap X_0 (g)) \times \GM}^{\GM}$.
\end{theorem}

\begin{proof}
Let $\varphi$ be in $\cL (X)$. A basic observation is that
when the inequality  $\ord_t f(\varphi)
< N \ord_t g (\varphi)$ holds,
$f(\varphi)$ and $(f+g^N)(\varphi)$ have same order
$\ord_t$ and same angular coefficient $\ac$. 
If $A$ is a subset of $\cL_n (X)$,
we  denote
by $A^+$, resp. $A^0$, the intersection of $A$ with the
set of arcs in  $\cL_n (X)$ such that $\ord_t f(\varphi) >
N\ord_t g(\varphi)$, resp. $\ord_t f(\varphi) = N\ord_t g
(\varphi)$. In this way one defines series
\begin{equation}
Z^+_f (T) = \sum_{n \geq 1} [\cX_n^+ (f)] \, \LL^{-n d} T^n
\end{equation}
and
\begin{equation}
Z^0_f (T) = \sum_{n \geq 1} [\cX_n^0 (f)] \, \LL^{-n d} T^n
\end{equation}
in $\cM_{X_0 (f) \times \GM}^{\GM} [[T]]$
and similarly series
$Z^+_{f + g^N} (T)$ and $Z^0_{f + g^N}(T)$
in $\cM_{X_0 (f + g^N) \times \GM}^{\GM} [[T]]$.
It follows from the previous remark that
\begin{equation}\label{last}
i_1^* Z_f (T) - i_2^* Z_{f + g^N} (T) = i_1^* (Z^+_f (T) + Z^0_f (T))
- i_2^* (Z^+_{f + g^N} (T) + Z^0_{f + g^N} (T)),
\end{equation}
where we extend $i^*_1$ and $i^*_2$ to series componentwise.

Let $N$ be a positive integer. For any integer $r$, we
denote by $\pi_N$ the morphism $X_0 (g) \times \GM^r \times
\GM \rightarrow X_0 (g) \times \GM^r\times \GM$ mapping
$(x,\mu, \lambda)$ to $(x,\mu , \lambda^N)$. Then we have
\begin{lem}\label{puissance}
Given a map $g:X\longrightarrow \AA^1_k$ and the induced
nearby cycles morphism $\cS_g$ from $\cM_{X_0 (g)\times
\GM^r}^{\GM^r}$ to $\cM_{X_0 (g)\times \GM^r\times
\GM}^{\GM^r\times \GM}$, then, for any positive integer
$N$, the following equality holds:
 \begin{equation}
\cS_{g^N}=\pi_{N!}\circ \cS_{g}.
 \end{equation}
\end{lem}
\begin{proof} 

Let $Z$ be a 
smooth variety with good $\GM^r$-action,
endowed with an equivariant  morphism
$p  : Z \rightarrow X$
and a monomial morphism ${\mathbf f} : Z \rightarrow \GM^r$ such that the morphism
$(p,  {\mathbf f}) : Z \rightarrow X \times \GM^r$ is proper,
and let $U = Z \setminus F$ be an open
dense subset of $Z$ which is
stable under the $\GM^r$-action. 
For a positive integer $\gamma$, with the notations of \ref{eqset},
let us consider the modified zeta function of $g^N$ on $U$ 
\begin{equation}
Z^{\gamma}_{g^N \circ p, U} (T) := 
\sum_{n\geq 1} \, 
[\cX^{\gamma n}_n (g^N \circ p, U)] \, \LL^{- \gamma n d} T^n
\end{equation} in $\cM_{X_0 (g)\times \GM^r \times \GM}^{\GM^r\times \GM}
[[T]]$.
Since 
$\cX^{\gamma n}_n (g^N \circ p, U)$, is empty unless $N$ divides $n$ and
\begin{equation}
[\cX^{\gamma mN}_{mN} (g^N \circ p, U)]
=
\pi_{N!} ([\cX^{\gamma m N}_m (g \circ p, U)]),
\end{equation}
we get that 
$Z^{\gamma}_{g^N \circ p, U} (T)$ is equal to  $\pi_{N!}(Z^{\gamma N}_{g \circ p,U}(T^N))$, the limit of
which, as $T$ goes to infinity, is equal, for $\gamma$ big
enough, to $\pi_{N !}(\cS_{g \circ p,U})$. The result follows from Theorem \ref{extequi}.
\end{proof}

\begin{lem}Assume
$N >  \gamma ((f), (g))$.

Then the series $i_2^* (Z^+_{f +
g^N} (T))$ lies in $\cM_{(X_0 (f) \cap X_0 (g)) \times
\GM}^{\GM} [[T]]_{\rm sr}$ and
\begin{equation}\label{calc1}
\lim_{T \mapsto \infty}  i_2^* (Z^+_{f + g^N} (T)) = - \cS_{g^N} ([X_0 (f)]).
\end{equation}
\end{lem}
\begin{proof}
Note that
 $\cX_{n}^+(f+g^N)$ is non empty only if $n$ is
a multiple of $N$ and 
that
\begin{equation}
\Bigl[ \cX_{mN}^+(f+g^N) \Bigr]= \pi_{N!} \Bigl(\Bigl[\Bigl \{\varphi \in \cL_{mN} (X) \Bigm
\vert \ord_t g (\varphi) = m,\; \ord_t f (\varphi) >
Nm\Bigr\}\Bigr]\Bigr),
\end{equation}
 the variety on the right hand side  being endowed with the morphism to $\GM$  induced by 
$\ac (g)$.
Summing up, we may write by (\ref{mzf}) and the proof of
Lemma
\ref{puissance}
\begin{equation}
 Z_{f+g^N}^+(T) =\pi_{N !}(Z_g (T^N) - Z^N_{g, X \setminus X_0 (f)} (T^N)).
\end{equation}
By Proposition \ref{formop} and its proof, for $N > \gamma
((f), (g))$, the series $Z^N_{g, X \setminus X_0 (f)} (T)$
lies in $\cM_{X_0 (g) \times \GM}^{\GM} [[T]]_{\rm sr}$ and
its $\lim_{T \mapsto \infty}$ is equal to $- \cS_{g, X
\setminus X_0 (f)}$.
The same holds for $Z^N_{g, X \setminus X_0 (f)} (T^N)$
and
the 
result follows since $\pi_{N !}
(\cS_g - \cS_{g, X \setminus X_0 (f)})$ is equal to
$\cS_{g^N} ([X_0 (f)])$ 
by Lemma \ref{puissance}.
\end{proof}

\subsection{} 
We fix an integer $N$ such that 
$N > \gamma ((f), (g))$ and   a log-resolution $h
: Y \rightarrow X$ of $(X, X_0 (f) \cup X_0 (g))$ such that $N > \gamma_h ((f), (g))$.
We keep the notations used in
\ref{resolutions} and \ref{resol}.  In particular, $N N_i (g)
> N_i (f)$ for $i \in C$. Note that the stratum
$E_I^{\circ}$ is contained in $(g \circ h)^{-1} (0)$ if and
only if $J = I \cap C$ is non empty.

Fix a non empty stratum  $E_I^{\circ}$ in $Y$.

We consider the cones $\Delta_I^+$ and $\Delta_I^0$ in
$\RR_{>0}^{I}$ defined respectively by
\begin{equation}
\sum_{i \in I} N_i (f)
x_i > N\sum_{j \in J} N_j
 (g) x_j
\end{equation}
and
\begin{equation} \sum_{i \in I} N_i
(f) x_i =N\sum_{j \in J} N_j
 (g) x_j.
\end{equation}
Note that when
$K = I \setminus C$ is empty, $I = J$ and
$N_i (f) < N N_i (g)$ for all $i$, hence the cones
$\Delta_I^+$ and $\Delta_I^0$ are both empty.

As in (\ref{trtr}), we have
\begin{equation}
i_1^* (Z^+_f  (T)+Z^0_f (T)) = \sum_{\sur{I\cap B \neq
\emptyset}{I \cap C\neq \emptyset}} [U_I] \, \Psi_I (T)
\end{equation}
with
\begin{equation}\Psi_I (T) = \sum_{{\bf k} \in (\Delta_I^+  \cup  \Delta_I^0)\cap \NN_{>
0}^I } \LL^{-\sum_{i \in I} \nu_i k_i} T^{\sum_{i \in I}
N_i (f) k_i}.
\end{equation}
Since
\begin{equation}
\lim_{T \mapsto \infty} \Psi_I (T) = \chi (\Delta_I^+  \cup \Delta_I^0) = 0,
\end{equation}
we deduce that
\begin{equation}\label{extra2}
\lim_{T \mapsto \infty} i_1^* (Z^+_f  (T)+Z^0_f (T)) = 0.
\end{equation}

\subsection{}\label{defcn}
We now want to compute the zeta function $Z^0_{f+g^N} (T)$. 
Fix
$\mathbf k$ in $\Delta_I^0 \cap \NN^I_{>0}$
and denote by $\phi$ 
the finite morphism from
$\AA^1_k$ to $\AA^I_k$
sending $u$ to  
$(u^{k_i})$.
We shall still denote by 
$\phi$ its restriction
as a group morphism from $\GM$ to $\GM^I$.
Taking the pullback along $\phi$ of the 
deformation to the normal cone to $E_I$ in $Y$,
$ p_I : CY_I \rightarrow  \AA_k^I$,
introduced in \ref{desdef}, one gets
a morphism
$p : CY_{\bf k} \rightarrow \AA^1_k$ having the following description.
The scheme $CY_{\bf k}$ may be identified with  $\Spec \cA_{\bf k}$ where 
\begin{equation}
\cA_{\bf k} := \sum_{{\bf n} \in \NN^I} \cO_{Y \times
\AA^1_k} \left(- \sum_{i \in I} n_i 
(E_i \times \AA^1_k)
\right ) u^{-
\sum_{i \in I} k_i n_i}
\end{equation}
is a subsheaf of $\cO_{Y \times \AA^1_k} [u^{-1}]$, and 
the natural inclusion $\cO_{Y \times \AA^1_k} \rightarrow
\cA_{\bf k}$ induces a morphism $\pi : CY_{\bf k}
\rightarrow Y \times \AA^1_k$ from which $p$ is derived. 
Via the same inclusion,
the functions $f\circ h$, $g^N\circ h$ and $(f + g^N) \circ
h$ are divisible by $ u^{\sum_{i \in
I}k_i N_i (f)}$
 in $\cA_{\bf k}$.
 We denote the quotients by $\widetilde
f_\mathbf{k}$, $\widetilde g_\mathbf{k}^N$ and $\widetilde
F_\mathbf{k}$, respectively.

We denote by $\widetilde E_i$ the 
pullback
 of the divisor
$E_i \times \AA^1_k$ by $\pi$, 
by
$D$ the divisor globally defined on $CY_\mathbf{k}$
by $u=0$, and by $CE_i$ the divisors $\widetilde E_i-k_iD$,
$i$ in $I$ (resp. $\widetilde E_i$,
$i$ not in $I$).
We denote by $CY_\mathbf{k}^\circ$ the
complement in $CY_\mathbf{k}$ of the union of the $CE_i$,
$i$ in $A$, and by $Y^\circ$ the complement in $Y$ of the
union of the $E_i$, $i$ in $A$. 
We denote by $F_I$ the
function $f_I + g_I^N :U_I \rightarrow \AA^1_k$.

\begin{lem} \label{lemcn}
The scheme $CY_{\bf k}$ is smooth, the morphism $\pi$
induces an isomorphism above $\AA^1_k \setminus \{0\}$, the
morphism $p$ is a smooth morphism and its fiber $p^{-1}
(0)$ may be naturally identified with the bundle
$\nu_{E_I}$. When restricted to $CY_{\bf k}^\circ$, the
fiber of $p$ above $0$ is naturally identified with 
$U_{I}$
and $\pi$ induces an isomorphism 
between 
 $CY_{\bf k}^\circ \setminus p^{-1}Ê(0)$ and $Y^\circ\times \AA^1_k
\setminus \{0\}$. 
The restriction of $\widetilde f_\mathbf{k}$
(resp. $\widetilde g_\mathbf{k}$, $\widetilde F_\mathbf{k}$) to the
fiber $U_I \subset p^{-1} (0)$ is equal to $f_I$ (resp.
$g_I$, $F_I$).
\end{lem}

\begin{proof}Since $CY_{\bf k}$ is covered by open subsets
of the form $\Spec \cO_U [y_i, u] / (z_i - u^{k_i} y_i)$
with $U$ open subset on which the divisors $E_i$ are
defined by equations $z_i = 0$, the smoothness of $CY_{\bf
k}$ is clear. The remaining properties are checked
directly.
\end{proof}

The $\GM^I$-action $\sigma_I$ on $CY_I$ induces
via $\phi$ a $\GM$-action on $CY_{\bf k}$ 
we denote by $\sigma$,  
leaving sections of $\cO_Y$ invariant and acting
on $u$
by
$\sigma (\lambda) : u \mapsto \lambda^{-1} u$.
 Note that in coordinate charts such  as in the proof of
 Lemma \ref{lemcn}, $\sigma$ leaves $z_i$ invariant and 
 $\sigma (\lambda)$
 maps $y_i$ to
 $\lambda^{k_i} y_i$. 
We have now two different $\GM$-actions
on $\cL_n(CY_\mathbf{k}^\circ)$: the one induced by the standard $\GM$-action
on arc spaces and the one induced by $\sigma$.
We denote by $\widetilde \sigma$ the action given by the composition of these
two (commuting) actions.

For $\varphi$ in $\cL_n (Y)$ with $\varphi (0)$ in $E_i$,
we set $\ord_{E_i} \varphi := \ord_t z_i (\varphi)$, for
$z_i$ any local equation of $E_i$ at $\varphi (0)$.

Let us denote by $\widetilde \cL_n(CY_\mathbf{k}^\circ)$ the
set of arcs $\varphi$ in $\cL_n(CY_\mathbf{k}^\circ)$ such
that $p(\varphi(t))=t$ (in particular $\varphi(0)$ is in $U_I$).
For such an arc $\varphi$, composition with $\pi$
sends
$\varphi$ to an arc in $\cL_n(Y\times\AA^1_k)$ which is the
graph of an arc in $\cL_n(Y)$ not contained in the union of
the divisors $E_i$, $i$ in $I$.
Note that $\widetilde \cL_n(CY_\mathbf{k}^\circ)$ is stable by $\widetilde \sigma$.

We will consider $\cX_{n,\mathbf{k}}$, the set of arcs
$\varphi$ in $\cL_n (Y)$ such that $\varphi(0)$ is in
$E_I^{\circ}$ and $\ord_{E_i} \varphi = k_i$ for $i \in I$.

\begin{lem}\label{fib}
Assume $n \geq k_i$ for $i$ in $I$.
The morphism
$\widetilde \pi
: \widetilde \cL_n(CY_\mathbf{k}^\circ) \rightarrow \cX_{n,\mathbf{k}}$
induced by the projection $CY_{\bf k}^\circ \rightarrow Y$ 
is an affine bundle 
with fiber $\AA^{\sum_I k_i}_k$.
Furthermore if  $\widetilde \cL_n(CY_\mathbf{k}^\circ)$
is endowed with the $\GM$-action induced by $\widetilde \sigma$ and
$\cX_{n,\mathbf{k}}$ with  the standard $\GM$-action,
$\widetilde \pi$ is $\GM$-equivariant and the action of $\GM$ on the affine bundle is affine.
If $n \geq \sum_I k_i N_i
(f)$, the
composed maps $\ac(f\circ h)(\widetilde \pi(\varphi))$ and
$\ac(g\circ h) (\widetilde \pi(\varphi))$ are equal
respectively to $f_{I}(\varphi(0))$ and
$g_{I}(\varphi(0))$, whereas
\begin{equation}
\ac ((f + g^N) \circ h) (\widetilde \pi (\varphi)) = \ac (\widetilde F_{\bf k} (\varphi)).
\end{equation}
Furthermore, when $F_I(\varphi(0))\neq 0$, hence $(\ord_t(f
+ g^N) \circ h)(\widetilde \pi (\varphi)) = \sum_I k_i N_i
(f)$, we have
\begin{equation}\ac ((f + g^N) \circ h) (\widetilde \pi
(\varphi))=F_I(\varphi(0)).
\end{equation}
\end{lem}

\begin{proof}Every point in $E_I^{\circ}$ is contained
in a open subset $U$ of $Y$ such that
the divisors $E_i$, $i \in I$
are defined by equations $z_i = 0$
in $U$ and such that there exists furthermore
$d - \vert I \vert$ functions $w_j$ on $U$
such that the family $(z_i, w_j)$ gives rise to
an \'etale morphism
$ U \rightarrow \AA^d_k$.  This morphism
induces an isomorphism
$\cL_n(U ) \simeq U \times_{\AA^d_k} \cL_n (\AA^d_k)$,  cf. Lemma 4.2 of  \cite{arcs}.
Adding further the coordinate $u$, gives an isomorphism
$\cL_n(U  \times \AA^1_k) \simeq (U \times \AA^1_k ) \times_{\AA^{d + 1}_k} \cL_n (\AA^{d + 1}_k)$.
The family $(y_i, w_j, u)$, with $z_i = y_i u^{k_i}$,
induces an \'etale morphism
$\pi^{-1} (U \times \AA^1_k) \rightarrow \AA^{d + 1}_k$,
hence an isomorphism
$\cL_n(\pi^{-1}(U  \times \AA^1_k)) \simeq (\pi^{-1} (U \times \AA^1_k)) \times_{\AA^{d + 1}_k} \cL_n (\AA^{d + 1}_k)$. Under
these isomorphisms $\widetilde \pi$ just corresponds to multiplicating
each $y_i$-component of an arc by
$t^{k_i}$.
Note in particular that in  that description the action of
$\widetilde \sigma (\lambda)$ on a component $y_i (t)$ is given by
$y_i (t) \mapsto \lambda^{k_i} y_i (\lambda t)$, hence
$\widetilde \pi$ is $\GM$-equivariant. 
The rest of the statement follows also directly from that description.
\end{proof}

We define $\cY_{n, {\bf k}}$ as the subset of
$\cX_{n,\mathbf{k}}$ consisting of those arcs $\varphi$ such that $\ord_t
((f + g^N)\circ h) (\varphi) = n$. The constructible set $\cY_{n,
{\bf k}}$ is stable by the usual $\GM$-action on
$\cL_n(Y)$ and the morphism $\ac(f+g^N)$ defines a class
$[\cY_{n, {\bf k}}]$ in $\cM_{(X_0 (f) \cap X_0 (g)) \times
\GM}^{\GM}$. 
By definition $\cY_{n, {\bf k}} = \emptyset$
if $n < \sum_I k_i N_i(f)$.

We then define $\widetilde \cY_{n, {\bf k}}$ as the 
preimage of $\cY_{n, {\bf k}}$ by the fibration $\widetilde\pi$ of
Lemma \ref{fib}. It consists of arcs $\varphi$ in
$\widetilde \cL_n(CY_{\mathbf{k}}^\circ)$ such that $\ord_t \widetilde
F_{\mathbf{k}}(\varphi)=n-\sum_I k_i N_i(f)$. 
We denote
by $[\widetilde \cY_{n, {\bf k}}]$ the class of $\widetilde \cY_{n,
{\bf k}}$ in $\cM_{(X_0 (f) \times X_0 (g)) \times
\GM}^{\GM}$, the morphism $\widetilde \cY_{n, {\bf k}}
\rightarrow \GM$ being $\ac (\widetilde F_{\bf k} )$ and the
$\GM$-action being induced by $\widetilde \sigma$.
We denote by $[U_I \setminus
(F_I^{-1} (0))]$ the class of $U_I \setminus (F_I^{-1}
(0))$ in $\cM_{(X_0 (f) \times X_0 (g)) \times \GM}^{\GM}$,
the $\GM$-action being the natural diagonal action of
weight ${\bf k}$ on $U_I \setminus (F_I^{-1} (0))$ and the
morphism to $\GM$ being the restriction of $F_I$. We also
consider the class $[\GM \times F_I^{-1} (0)]$ of $\GM
\times F_I^{-1} (0)$ in $\cM_{(X_0 (f) \times X_0 (g))
\times \GM}^{\GM}$, the $\GM$-action on the second factor
being the diagonal one and the morphism to $\GM$ being the
first projection.

\begin{lem}\label{cal}The following equalities
hold  in $\cM_{(X_0 (f) \times X_0 (g)) \times \GM}^{\GM}$
\begin{enumerate}
\item[(1)] $[\widetilde \cY_{n, {\bf k}}] =\LL^{nd} \, [U_I \setminus (F_I^{-1} (0))]$,
if $n = \sum_I k_i N_i (f)$,
\item[(2)]$[\widetilde \cY_{n, {\bf k}}] = \LL^{nd -m} \, [\GM \times F_I^{-1} (0)]$,
if $n - \sum_I k_i N_i (f) = m> 0$.
\end{enumerate}
\end{lem}

\begin{proof}If $n = \sum_I k_i N_i (f)$,
$\widetilde \cY_{n, {\bf k}}$ is the set of arcs $\varphi (t)$
in $\cL_n (CY_{\bf k}^\circ)$ such that $\varphi (0)$ lies
in $U_I \setminus (F_I^{-1} (0))$  and $u (\varphi (t)) = t$, and $(1)$ follows.

If $n - \sum_I k_i N_i (f) = m> 0$, $\widetilde \cY_{n, {\bf
k}}$ is the set of arcs $\varphi $ in $\cL_n (CY_{\bf
k}^\circ)$ such that $\ord_t (\widetilde F_{\bf k})( \varphi) = m$ and $u (\varphi (t)) = t$. Now let us observe that the
morphism $(\widetilde F_{\bf k}, u) : CY_{\bf k}^\circ
\rightarrow \AA^2_k$ is smooth on a neighborhood of $U_I$
in $CY_{\bf k}^\circ$,  
since
$u$ is a smooth function on
$CY_{\bf k}^\circ$ and 
the restriction of $\widetilde F_{\bf
k}$ to the divisor $u=0$, identified with $U_I$, is
$F_I=f_I+g_I^N$ which is a smooth function on $U_I$. The fact that 
$F_I=f_I+g_I^N$ is a smooth function on $U_I$ is checked locally as follows:
for $i$ in $I \setminus C$ (recall $I \setminus C$ is non empty), and with local coordinates as above,
\begin{equation}
y_i \frac{\partial}{\partial y_i}Ê(f_I + g_I^N) = y_i \frac{\partial f_I}{\partial y_i} = N_i (f) f_I
\end{equation}
does not vanish on $U_I$.
\end{proof}

\begin{lem}\label{cvc}We have
\begin{equation}
i_2^* [\cX_n^0 (f + g^N)]= \sum_{\sur{I\cap
C=J\neq \emptyset}{I\setminus C=K\neq \emptyset}}
\sum_{{\bf k} \in \Delta_I^0 \cap \NN_{> 0}^I }
[\cY_{n, {\bf k}}] \, \LL^{- \sum_{i \in I} (\nu_i - 1) k_i}.
\end{equation}
\end{lem}

\begin{proof}This is a standard application of the change of
variable formula, or more precisely of Lemma 3.4 in
\cite{arcs}. The proof is completely similar to the proof
of Theorem 2.4 of \cite{lef}. (Recall that $\Delta_I^0$ is
empty if $K$ is empty.)
\end{proof}
It follows from Lemma \ref{cvc} and Lemma \ref{fib} that
\begin{equation}
i_2^* Z^0_{f + g^N} (T) = \sum_{n >0}
\sum_{\sur{I\cap
C=J\neq \emptyset}{I\setminus C=K\neq \emptyset}}
\sum_{{\bf k} \in \Delta_I^0 \cap \NN_{> 0}^I}
[\widetilde \cY_{n, {\bf k}}] \, \LL^{- \sum_{i \in I} \nu_i k_i}\LL^{-nd} \,T^n.
\end{equation}
Using Lemma \ref{cal}, we
deduce
\begin{equation}\label{yui}
i_2^* Z^0_{f + g^N} (T) = \sum_{\sur{I\cap
C=J\neq \emptyset}{I\setminus C=K\neq \emptyset}}
\Bigl(
[U_I \setminus (F_I^{-1} (0))] +
 [\GM \times F_I^{-1} (0))]\frac{\LL^{-1} T}{1 - \LL^{-1} T}
\Bigr)
\Phi_I (T),
\end{equation}
with
\begin{equation}
\Phi_I (T) = \sum_{{\bf k} \in \Delta_I^0 \cap \NN_{> 0}^I}
\LL^{- \sum_{i \in I} \nu_i  k_i} T^{\sum_{i \in I} N_i (f) k_i}.
\end{equation}
Since $\Phi_I (T)$ lies in
$\ZZ [\LL, \LL^{-1}] [[T]]_{\rm sr}$ and
\begin{equation}\label{limphi}
\lim_{T \mapsto \infty} \Phi_I (T) = \chi (\Delta_I^0) = (-1)^{\vert I \vert - 1},
\end{equation}
$i_2^* Z^0_{f + g^N} (T) $ lies in $\cM_{(X_0(f) \cap X_0
(g)) \times \GM}^{\GM} [[T]]_{\rm sr}$. Using Theorem
\ref{reee}, we deduce from (\ref{yui}) and  (\ref{limphi})
that
\begin{equation}\label{extra}
\lim_{T \mapsto \infty}i_2^* Z^0_{f + g^N} (T) = \Psi_{\Sigma} (\cS_{g^N} (\cS_f)).
\end{equation}

We deduce from  (\ref{last}), (\ref{calc1}), (\ref{extra}) and (\ref{extra2})
that
\begin{equation}\label{fini}
i_1^* \cS_f - i_2^* \cS_{f + g^N} = \Psi_{\Sigma} (\cS_{g^N} (\cS_f)) - \cS_{g^N} ([X_0 (f)]).
\end{equation}

By Proposition \ref{ass}, $\cS_{g^N} ([X_0 (f)]) = \Psi_{\Sigma} (\cS_{g^N} ([\GM \times X_0 (f)]))$,
hence the statement of the Theorem
follows from (\ref{fini}), since
$i_1^* \cS^{\phi}_f - i_2^* \cS^{\phi}_{f + g^N} = (-1)^{d - 1}(i_1^* \cS_f - i_2^* \cS_{f + g^N})$.
\end{proof}

If $f$ is a function on the smooth variety $X$ of pure dimension $d$ and
$x$
is a closed point of $X_0 (f)$, we write $\cS_{f, x}$ for $i_x^*\cS_f$,
and $\cS^{\phi}_{f, x}$ for $i_x^*\cS^{\phi}_f$,
where $i_x$ stands for the inclusion of $x$ in $X_0 (f)$.
Note that $\cS^{\phi}_{f, x} = (-1)^{d - 1} (\cS_{f, x} -[ \GM \times \{x\}])$.

Theorem \ref{MT} has the following local corollary:

\begin{cor}\label{mc}
Let $X$ be a smooth variety of pure dimension $d$, and $f$
and $g$ be two functions  from $X$ to $\AA^1_k$. Let $x$ be
a closed point of $X_0 (f) \cap X_0 (g)$. For every $N >
\gamma_x ((f), (g))$, the equality 
\begin{equation}\label{mtc}
 \cS_{f, x}^\phi - \cS_{f+g^N, x}^\phi = \Psi_{\Sigma}
(\cS_{g^N, x}(\cS_f^\phi))
\end{equation}holds in
$\cM_{\GM}^{\GM}$.
\end{cor}

\begin{proof}
The only point to be checked is that $\gamma ((f), (g))$
may be replaced by the local invariant $\gamma_x ((f),
(g))$, which is clear from the proof of Theorem
\ref{MT}.
\end{proof}

\subsection{}Let us now explain how to deduce from
Theorem \ref{MT} the motivic Thom-Sebastiani Theorem
of \cite{ts}, \cite{looi} and \cite{barc}.

Let $X$ and $Y$ be two varieties over $k$. For $r$ and $s$ in $\NN$,
cartesian product gives rise to an external product
\begin{equation}
\boxtimes : \cM_{X \times \GM^r}^{\GM^r} \times
\cM_{Y \times \GM^s}^{\GM^s}
\longrightarrow
\cM_{X \times Y \times \GM^{r + s}}^{\GM^{r + s}}
\end{equation}
(not to be confused with the one in \ref{noto}).

Let $X_1$ and $X_2$ be smooth varieties
and consider
functions $f_1 : X_1 \rightarrow \AA^1_k$ and $f_2 : X_2 \rightarrow \AA^1_k$.
We set $X_0 = f_1^{-1} (0) \times f_2^{-1} (0)$
and, for any
$Y \subset X_1 \times X_2$ containing $X_0$ 
we denote by $i$ the inclusion
of $X_0 \times \GM$ in $Y \times \GM$.

\begin{theorem}\label{thomseb}Let $X_1$ and $X_2$ be smooth varieties
of pure dimension $d_1$ and $d_2$ 
and consider
functions $f_1 : X_1 \rightarrow \AA^1_k$ and $f_2 : X_2 \rightarrow \AA^1_k$.
Denote by
$f _1\oplus f_2$ the function on
$X_1 \times X_2$ sending $(x_1, x_2)$ to $f_1 (x_1) + f_2 (x_2)$.
Then
\begin{equation}\label{mur}
i^* \cS^{\phi}_{f_1 \oplus f_2} = 
 \Psi_{\Sigma} (\cS^{\phi}_{f_1} \boxtimes  \cS^{\phi}_{f_2})
\end{equation}
in $\cM_{X_0 \times \GM}^{\GM}$.
\end{theorem}

\begin{proof}We set $X = X_1 \times X_2$
and we denote by $f$ and $g$
the functions on $X$ induced by
$f_1$ and $f_2$, respectively.
In particular $f_1 \oplus f_2 = f + g$.
If $Y_1\rightarrow X_1$ is a log-resolution of
$(X_1, f_1^{-1} (0))$ and
 $Y_2\rightarrow X_2$ is a log-resolution of
$(X_2, f_2^{-1} (0))$,
$h : Y := Y_1 \times Y_2 \rightarrow X$
is a log-resolution of
$(X, f^{-1} (0) \cup g^{-1} (0))$.
Using such a  log-resolution it is easily checked that
$\gamma ((f), ( g)) = 0$.
By (\ref{fini}),
\begin{equation}\label{encoreune}
i^* \cS_{f} - i^* \cS_{f + g} = \Psi_{\Sigma}( \cS_{g} (\cS_f)) - i^* \cS_{g} ([X_0 (f)]).
\end{equation}
Using the log-resolution $h$ one checks
that
$i^* \cS_{f} = \cS_{f_1} \boxtimes [f_2^{-1} (0)]$,
$ \cS_{g} (\cS_f) = \cS_{f_1} \boxtimes \cS_{f_2}$
and
$i^* \cS_{g} ([X_0 (f)]) = [f_1^{-1} (0) ]
\boxtimes \cS_{f_2}$.
Hence (\ref{encoreune}) may be rewritten as
\begin{equation}
 \Psi_{\Sigma} (\cS_{f_1} \boxtimes  \cS_{f_2})
 =
 \cS_{f_1} \boxtimes [f_2^{-1} (0)]
+
[f_1^{-1} (0)] \boxtimes \cS_{f_2}
- i^* \cS_{f_1 \oplus f_2}.
\end{equation}
Since $\cS_{f_1} \boxtimes [f_2^{-1} (0)] = \Psi_{\Sigma} (\cS_{f_1} \boxtimes [f_2^{-1} (0) \times \GM])$
and 
$[f_1^{-1} (0)] \boxtimes \cS_{f_2}
= \Psi_{\Sigma} ([f_1^{-1} (0) \times \GM] \boxtimes \cS_{f_2})$
(cf. the proof of the statement concerning the unit element in Proposition \ref{ass}),
 (\ref{mur}) directly follows, by definition of $\cS^{\phi}$. 
 \end{proof}

\section{Spectrum and the Steenbrink conjecture}\label{ssc}

\subsection{}We now assume $k = \CC$.
We denote by $\HS$ the abelian category of Hodge
structures and by $K_0 (\HS)$ the corresponding Grothendieck ring
(see, eg, \cite{barc} for definitions).
Note that any mixed Hodge structure has
a canonical class in $K_0 (\HS)$.
Recall there is a canonical morphism
\begin{equation}
\chi_h : \cM_{\CC} \longrightarrow
K_0 (\HS),
\end{equation}
which assigns to the class a variety $X$ the element
$\sum_i (-1)^i [H^i_c(X,\QQ)]$ in $K_0(\HS)$,
where $[H^i_c(X,\QQ)]$ stands for  the class of the mixed Hodge structure
on $H^i_c(X,\QQ)$.
Let us denote by
$\HS^{\rm mon}$ the abelian category of
Hodge
structures endowed  with an automorphism of finite order
and by
$K_0 (\HS^{\rm mon})$
the corresponding Grothendieck ring.
Let us consider the ring morphism
\begin{equation}
\chi_h : \cM_{\GM}^{\GM} \longrightarrow
K_0 (\HS^{\rm mon})
\end{equation}
deduced from (\ref{tro}) via (\ref{comp})
and composition
with $K_0 ({\rm MHM}^{\rm mon}_{\Spec \CC}) \rightarrow 
K_0 (\HS^{\rm mon})$.
It is described as follows:  if  $[X]$ is the class
of $f : X \rightarrow \GM$ in $\cM_{\GM}^{\GM}$
with $X$ connected,
since $f$ is monomial with respect to the $\GM$-action,
$f$ is a locally trivial fibration for the complex topology.
Furthermore, if the weight is, say, $n$,
$x \mapsto \exp (2 \pi i t /n ) x$ is a geometric
monodromy of finite order along the origin.
It follows that $X_1$, the fiber of $f$ at $1$, is endowed with 
an automorphism of finite order $T_f$, and we have
\begin{equation}\label{9696}
\chi_h ([f : X \rightarrow \GM]) = (\sum_i (-1)^i [H^i_c(X_1,\QQ)], T_f).
\end{equation}

There is a natural linear map, called the Hodge spectrum,
\begin{equation}
\hsp :
K_0(\HS^{\rm mon }) \longrightarrow \ZZ[\QQ],
\end{equation}
such that
\begin{equation}\hsp ([H]) := \sum_{\alpha \in \QQ \cap[0,1)} t^\alpha (\sum_{p,q
\in
\ZZ} \dim(H^{p,q}_{\alpha})t^p),
\end{equation}
for any Hodge structure $H$ with an
automorphism of finite order, where $H^{p,q}_\alpha$ is the
eigenspace of $H^{p,q}$ with respect to the eigenvalue
$\exp (2\pi i\alpha)$.
We identify here  $\ZZ[\QQ]$ with $ \cup_{n
\geq 1} \ZZ[t^{1/n},t^{-1/n}]$.

We shall consider the composite morphism
\begin{equation}
\Sp := (\hsp \circ \chi_h) :
\cM_{\GM}^{\GM} \longrightarrow
\ZZ[\QQ].
\end{equation}
Note that $\Sp$ is a ring morphism for the convolution
product $*$ on $\cM_{\GM}^{\GM}$,
 by Lemma \ref{wfermat}.

Denoting by
$\HS^{2-{\rm mon}}$ the abelian category of
Hodge structures
endowed  with two commuting automorphisms of finite order
and by $K_0 (\HS^{2-{\rm mon}})$
the corresponding Grothendieck ring, one
deduce from (\ref{tro}) via (\ref{comp}) a ring
morphism
\begin{equation}
\chi_h : \cM_{\GM^2}^{\GM^2} \longrightarrow
K_0 (\HS^{2-{\rm mon}})
\end{equation}
having a description similar to  (\ref{9696}).

Also we can define a Hodge spectrum on $K_0(\HS^{2-{\rm mon}})$
as follows.
Denote by $\pi : [0, 1)\cap \QQ \rightarrow \QQ/\ZZ$ the restriction of
the projection $\QQ \rightarrow \QQ/\ZZ$ and by
$s : \QQ/\ZZ \rightarrow [0, 1) \cap \QQ$ its inverse.
The bijection $\QQ / \ZZ \times \ZZ \rightarrow \QQ$ sending
$(a, b) $ to $s (a) + b$ induces an isomorphism
of abelian groups between $\ZZ [\QQ / \ZZ \times \ZZ]$ and  $\ZZ[\QQ]$.
We define the spectrum
\begin{equation}
\hsp :
K_0(\HS^{2-{\rm mon}}) \longrightarrow \ZZ [(\QQ / \ZZ)^2 \times \ZZ]
\end{equation}
by
\begin{equation}\label{2sp}\hsp ([H]) = \sum_{\alpha\in
\mathbf{Q}\cap [0,1)}
\sum_{\beta\in \mathbf{Q}\cap
[0,1)}
 \sum_{p,q\in
\mathbf{Z}}
\Bigl(\dim H^{p,q}_{\alpha,\beta}\Bigr)
t^{\pi (\alpha)}
u^{\pi (\beta)}
 v^p,
 \end{equation}
with $H^{p,q}_{\alpha, \beta}$ the
eigenspace of $H^{p,q}$ with respect to the eigenvalue
$\exp (2\pi i\alpha)$ for the the first automorphism
and
$\exp (2\pi i\beta)$ for the the second  automorphism.
We shall denote by
$\Sp$ the morphism of abelian groups
\begin{equation}
\Sp := (\hsp \circ \chi_h) :
\cM_{\GM^2}^{\GM^2} \longrightarrow
\ZZ [(\QQ / \ZZ)^2 \times \ZZ].
\end{equation}

We denote by $\delta$
the morphism of abelian groups
\begin{equation}
\ZZ [(\QQ / \ZZ)^2 \times \ZZ]
\longrightarrow \ZZ [\QQ]
\end{equation}
sending
$t^a u^b v^c$
to
$t^{s (a) + s (b) + c}$.

Let $A$ be an element of
$\cM_{\GM^2}^{\GM^2}$. The relation between the spectrum of
$A$ and the spectrum of
$\Psi_{\Sigma} (A)$ is given by  the following
proposition.

\begin{prop}\label{compspec}
Let $A$ be an element of $\cM_{\GM^2}^{\GM^2}$. We have
\begin{equation}
\Sp (\Psi_{\Sigma} (A)) = \delta (\Sp (A)).
\end{equation}
\end{prop}
\begin{proof}
Let $A$ be a smooth variety with a good $\GM^2$-action
and with a morphism to $\GM^2$
which is diagonally monomial of weight $(n,n)$, $n$ in
$\NN_{>0}$. Let us denote by $A_1$ the fiber of $A$ above $(1,1)$.
By (\ref{prefermat}) and using the notation therein, we
have
\begin{equation}\label{ddd}
\Psi_{\Sigma} 
([A]) = - [F_1^n \times^{\mu_n \times \mu_n} A_1]
+
[F_0^n \times^{\mu_n \times \mu_n} A_1].
\end{equation}
The result follows from 
the following  well-known computation 
of the cohomology of Fermat varieties 
(cf. \cite{sk} and Lemma 7.1 in \cite{looi}).
\end{proof}

\begin{lem}\label{wfermat}
Let  $(\alpha ,\beta)$ be in $(\QQ/\ZZ)^2$. For every common 
denominator $n$ of $\alpha$ and $\beta$, the Hodge type of the 
eigenspaces
$H^i_c(F^n_1, \CC) (\alpha, \beta)$ and
$H^i_c(F^n_0, \CC) (\alpha, \beta)$
of $\mu_n\times\mu_n$ in $H^i_c(F^n_1, \CC)$ and
$H^i_c(F^n_0, \CC)$ respectively
with
character $(\alpha ,\beta)\in (n^{-1}\ZZ/\ZZ)^2$ is independent of $n$
and is computed as follows:
\begin{enumerate}
\item[(1)] 
$H^1_c(F^n_1, \CC) (\alpha, \beta)$ is of rank 1
for $(\alpha ,\beta)\not=(0,0)$ and  
and of rank 2 for
$(\alpha ,\beta) =(0,0)$ ;
$H^1_c(F^n_1, \CC) (\alpha, \beta)$ is of Hodge type 
$(0,1)$ if $\alpha\not= 0\not=\beta$  and 
$0< s (\alpha) +s (\beta) <1$,
$(1,0)$ if $1< s(\alpha) +s (\beta) <2$ and 
$(0,0)$ otherwise, that is if $\alpha=0$ or $\beta=0$ or
$\alpha +\beta=0$ ; 
$H^2_c(F^n_1) (0, 0)$ is of rank 1 and
Hodge type $(1, 1)$ ; 
all other cohomology groups are zero. 
\item[(2)] 
$H^1_c(F^n_0, \CC) (\alpha, - \alpha)$,
resp. 
$H^2_c(F^n_0, \CC) (\alpha, - \alpha)$, is of rank 1
and Hodge type $(0, 0)$, resp. 
$(1, 1)$, for any $\alpha$ in $\QQ / \ZZ$, 
and all other cohomology groups are zero. \qed
\end{enumerate}
\end{lem}

We shall also need the following
obvious statement:

\begin{lem}\label{supertriv}For $N \geq 1$, consider
the morphism
$\pi_N : \GM^2 \rightarrow \GM^2$ given by $(a, b) \mapsto
(a, b^N)$.
For every $A$ in
$\cM_{\GM^2}^{\GM^2}$,
\begin{equation}
\Sp (\pi_{N !} (A)) = \frac{1 - u}{1- u^{\frac{1}{N}}}\Sp  (A) (t, u^{\frac{1}{N}}, v).
\end{equation}
\end{lem}

\subsection{}\label{3.1.3}
Let $X$ be a smooth complex algebraic variety of dimension $d$
and let $f$ be a function $X \rightarrow \AA^1$.
Fix a closed point $x$ of $X$ at which $f$ vanishes.
Denote by $F_x$ the Milnor fiber of $f$ at $x$.
The cohomology groups $H^i (F_x,\QQ)$
carry a natural mixed Hodge structure
(\cite{St1}, \cite{Va}, \cite{Sa1}, \cite{Sa2}),
which
is compatible with the semi-simplification
of the monodromy operator $T_{f, x}$.  Hence we can define the Hodge
characteristic $\chi_h(F_x)$ of $F_x$ in
$K_0(\HS^{\rm mon})$.
The following statement follows from
\cite{motivic} and \cite{barc} (it is also a consequence of  Proposition \ref{com}):

\begin{theorem}\label{3.5.5}
Assuming the previous notations,
the following equality holds in
$K_0({\HS^{\rm mon}})$:
\begin{equation}
\chi_h(F_x) = \chi_h(\cS_{f, x}).
\end{equation}
\end{theorem}

In particular, if we
define the Hodge spectrum of $f$ at $x$
as
\begin{equation}\label{defhod}
\Sp (f, x)  := (-1)^{d- 1}\hsp (\chi_h(F_x)-1),
\end{equation}
it follows
from Theorem \ref{3.5.5}
that
\begin{equation}\label{analogy}
\Sp (f,x)  = \Sp (\cS^\phi_{f, x}).
\end{equation}

Now if $g : X \rightarrow \AA^1$
is another function vanishing at $x$,
we  shall set,
by analogy with (\ref{analogy}),
\begin{equation}
\Sp (f, g, x) := \Sp (\cS_{g, x}(\cS^\phi_{f})).
\end{equation}

Let us  denote by $\delta_N$
the morphism of abelian groups
$\ZZ [(\QQ / \ZZ)^2 \times \ZZ]
\rightarrow \ZZ [\QQ]$
sending
$t^a u^b v^c$
to
$t^{s (a) + s (b)/N + c}$.

\begin{prop}\label{hodge}For every positive integer $N$, the spectrum of
$\Psi_{\Sigma} (\cS_{g^N,x} (\cS_f^\phi))$ is equal to
\begin{equation}
\Sp(\Psi_{\Sigma}(\cS_{g^N,x} (\cS_f^\phi))) = \frac{1-t}{1-t^{\frac{1}{N}}}
\delta_N (\Sp (f,g,x)).
\end{equation}
\end{prop}

\begin{proof}Follows directly from Proposition \ref{compspec}
and Lemma \ref{supertriv}.
\end{proof}

Hence, we deduce immediately the following statement from
Corollary \ref{mc}.

\begin{theorem}Let $X$ be a smooth variety of pure dimension $d$,
and $f$ and $g$ be two functions  from $X$
to $\AA^1$.
Let $x$ be a closed point of $X_0 (f) \cap X_0 (g)$.
Then, for $N > \gamma_x ((f), (g))$,
\begin{equation}
\Sp (f,x) -
\Sp (f+g^N,x)=\frac{1-t}{1-t^{\frac{1}{N}}}\;
\delta_N (\Sp (f,g,x)).
\end{equation}
\end{theorem}

\subsection{Application to Steenbrink's conjecture}
Let us assume now
that the function $g$ vanishes on all local components at $x$
of the singular locus of $f$ but a finite number of
locally irreducible curves
$\Gamma_\ell$, $1\leq \ell \leq r$. We denote
by $e_{\ell}$ the order of
$g$ on $\Gamma_{\ell}$.

As in the introduction,
along the complement $\Gamma_{\ell}^\circ$
to $\{x\}$
in $\Gamma_{\ell}$, we may view
$f$ as a family
of isolated hypersurface singularities
parametrized by $\Gamma_{\ell}^\circ$.
We denote by $\alpha_{\ell, j}$ the exponents of
that isolated hypersurface singularity
and we note that there are two commuting
monodromy actions on 
the cohomology of
its Milnor fiber: the first one denoted by
$T_f$
is induced transversally
by the monodromy action of $f$ and the second one denoted by
$T_{\tau}$ is the monodromy
around $x$ in $\Gamma_{\ell}^\circ$.
Since
the semi-simplifications of $T_f$ and $T_{\tau}$ can be simultaneously diagonalized,
we may define
rational numbers
$\beta_{\ell, j}$ in $[0, 1)$
so that each $\exp (2 \pi i \beta_{\ell, j})$
is the eigenvalue of the semi-simplification of $T_{\tau}$
on the eigenspace of the semi-simplification of $T_{f}$
associated to $\alpha_{\ell, j}$.

We may now deduce from Theorem \ref{MT}, the following statement,
first proved by M. Saito in \cite{Sa3}, and later given another proof by
A. N\'emethi and J. Steenbrink in \cite{ns}.

\begin{theorem}\label{nstee}For $N >
\gamma_x ((f), (g))$, we have
\begin{equation}\label{uuu}
\Sp (f + g^N, x)
-
\Sp (f, x) = \sum_{\ell, j}
t^{\alpha_{\ell, j} + (\beta_{\ell, j}/e_{\ell} N)}
\frac{1- t}{1 -t^{1/e_{\ell} N}}.
\end{equation}
\end{theorem}


\begin{proof}
For every $\ell$, we set $\cS^{\phi}_{f, \ell} := i_{\ell}^{\ast} (\cS^{\phi}_f)$,
with $i_{\ell}$
the inclusion of $\Gamma_{\ell}^{\circ}$ in $X_0 (f)$.
Since $\cS^{\phi}_f - \sum_{\ell}  i_{\ell !} (\cS^{\phi}_{f, \ell})$
has support in $X_0 (g) \times \GM$,
\begin{equation}\label{ouf}
\cS_{g^N, x} (\cS^{\phi}_f) =  \sum_{\ell}  \cS_{g^N, x} (i_{\ell !} (\cS^{\phi}_{f, \ell})).
\end{equation}
Now consider the normalization $n_{\ell} :
\widetilde{\Gamma_{\ell}} \rightarrow \Gamma_{\ell}$.
Let us choose a 
uniformizing
parameter $\tau_{\ell}$ at the preimage $x_{\ell}$ of $x$ in 
$\widetilde{\Gamma_{\ell}}$. We may write
$g \circ n_{\ell} =  \eta  \tau_{\ell}^{e_{\ell}}$, with $\eta $ a local unit.
We have
\begin{equation}\label{qqq}
\cS_{g^N, x} (i_{\ell !} (\cS^{\phi}_{f, \ell})) = \cS_{g^N \vert \Gamma_{\ell}, x}
(\cS^{\phi}_{f, \ell}) = \cS_{\eta  \tau_{\ell}^{e_{\ell}N}, x_{\ell}} (\cS^{\phi}_{f, \ell}),
\end{equation}
where in the last term we view $\cS^{\phi}_{f, \ell}$ as lying in
$\cM_{\widetilde{\Gamma_{\ell}} \times \GM}^{\GM}$.
By Proposition  \ref{com},
\begin{equation}\label{yyy}
\Sp (\cS_{\eta \tau_{\ell}^{e_{\ell}N}, x_{\ell}} (\cS^{\phi}_{f, \ell})) =  \Sp (\cS_{\tau_{\ell}^{e_{\ell}N}, x_{\ell}} (\cS^{\phi}_{f, \ell}))
\end{equation}
and
\begin{equation}\label{ppp}
\Sp (\cS_{\tau_{\ell}, x_{\ell}} (\cS^{\phi}_{f, \ell})) = - \sum_j t^{\pi (\alpha_{\ell, j})} u^{\pi (\beta_{\ell, j})}
v^{[\alpha_{\ell, j}]},
\end{equation}
where $[ \alpha ]$ denotes the integer part of $\alpha$.
Indeed, 
note that if $H$ is the mixed Hodge module  corresponding
to a 
variation of mixed Hodge
structure on a neighborhood
of $x_{\ell}$, the fiber at $x_{\ell}$ of $\psi_{\tau_{\ell}} (H)$
is nothing but the generic fiber of the variation endowed with the monodromy around
$x_{\ell}$. The sign in (\ref{ppp})  results from the fact that 
the numbers
$\alpha_{\ell, j}$ occuring in its right hand side
are  the exponents of
an isolated hypersurface singularity in an ambient space of dimension $d - 1$ and not $d$.
The result follows now from (\ref{ouf}),
(\ref{qqq}), (\ref{yyy})  and
(\ref{ppp}) by
plugging together
Corollary \ref{mc} and 
Proposition \ref{hodge}.
\end{proof}

\bibliographystyle{amsplain}

\begin{thebibliography}{}





\bibitem{bm}E. Bierstone, P. Milman,
\textit{Canonical resolution of singularities in characteristic zero by
blowing up the maximal strata of a local invariant},
Invent. Math. \textbf{128} (1997), 207--302.

\bibitem{bittner}F. Bittner,
\textit{The universal Euler characteristic for varieties of
characteristic zero}, Compositio Math. \textbf{140} (2004), 1011--1032.



\bibitem{bittner2}F. Bittner,
\textit{On motivic zeta functions and the motivic 
nearby
fiber}, Math. Z. \textbf{249} (2005), 63--83.


\bibitem{D89}
J. Denef, \textit{On the degree of Igusa's local zeta
function}, Amer. J. Math. \textbf{109} (1987), 991--1008.

\bibitem{motivic}
J. Denef, F. Loeser, \textit{Motivic Igusa zeta functions},
J. Algebraic Geom. \textbf{7} (1998), 505--537.



\bibitem{arcs}J. Denef, F. Loeser,
\textit{Germs of arcs on singular algebraic varieties and motivic integration},
Invent. Math. \textbf{135} (1999), 201--232.

\bibitem{ts}J. Denef, F. Loeser,
\textit{Motivic exponential integrals and a motivic
Thom-Sebastiani Theorem}, Duke Math. J. \textbf{99} (1999),
285--309.


\bibitem{barc}
J. Denef, F. Loeser, \textit{Geometry on arc spaces of
algebraic varieties}, Proceedings of 3rd European Congress
of Mathematics, Barcelona 2000, Progress in Mathematics
\textbf{201} (2001), 327--348, Birkha{\"u}ser.

\bibitem{lef}J. Denef, F. Loeser,
\textit{Lefschetz numbers of iterates of the monodromy and
truncated arcs}, Topology \textbf{41} (2002), 1031--1040.



\bibitem{eh}
S. Encinas, H. Hauser,
\textit{Strong resolution of singularities in characteristic zero},
Comment. Math. Helv.
\textbf{77} (2002), 821--845.


\bibitem{ev}
S. Encinas, O. Villamayor,
\textit{Good points and constructive resolution of singularities},
Acta Math.  \textbf{181} (1998), 109--158.



\bibitem{ful}
W. Fulton, \textit{Intersection theory},
Ergebnisse der Mathematik und ihrer Grenzgebiete,
Springer-Verlag, Berlin, 1984.


\bibitem{guibert}G. Guibert,
\textit{Espaces d'arcs et invariants d'Alexander}, Comment.
Math. Helv. \textbf{77} (2002), 783--820.

\bibitem{iomdin}I.N. Iomdin,
\textit{Complex surfaces with a one-dimensional set of singularities}
(Russian),
Sibirsk. Mat. {\v Z.} \textbf{15} (1974), 1061--1082, 1181,
English translation: Siberian Math. J.  \textbf{15} (1974), no. 5,
748--762 (1975).



\bibitem{lmb}
G. Laumon,  L. Moret-Bailly,
\textit{Champs alg\'ebriques},
Ergebnisse der Mathematik und ihrer Grenzgebiete,
Springer-Verlag, Berlin, 2000.

\bibitem{looi}
E. Looijenga, \textit{Motivic Measures}, Ast\'erisque
\textbf{276} (2002), 267--297, S{\'e}minaire Bourbaki,
expos{\'e} 874.

\bibitem{ns}
A. N{\'e}methi, J. Steenbrink, \textit{Spectral pairs,
mixed Hodge modules, and series of plane curve
singularities}, New York J. Math. \textbf{1} (1994/95),
149--177.



\bibitem{Sa1}M. Saito, \textit{Modules de Hodge polarisables}, Publ. Res.
Inst. Math. Sci. \textbf{24} (1988), 849--995.




\bibitem{Sa1b}M. Saito, \textit{Duality for vanishing cycle functors}, Publ.
Res. Inst. Math. Sci. \textbf{25} (1989), 889--921.



\bibitem{Sa2}M. Saito, \textit{Mixed Hodge modules}, Publ.
Res. Inst. Math. Sci. \textbf{26} (1990), 221--333.



\bibitem{Sa3}M. Saito,
\textit{On Steenbrink's conjecture}, Math. Ann.
\textbf{289} (1991), 703--716.


\bibitem{sk}T. Shioda, T. Katsura,
\textit{On Fermat varieties}, T{\^o}hoku Math. J.
\textbf{31} (1979), 97--115.


\bibitem{Siersma}D. Siersma,
\textit{The monodromy of a series of hypersurface
singularities}, Comment. Math. Helv. \textbf{65} (1990),
181--197.


\bibitem{St1}J. Steenbrink,
\textit{Mixed Hodge structures on the vanishing
cohomology}, in Real and Complex Singularities, Sijthoff
and Noordhoff, Alphen aan den Rijn, 1977, 525--563.








\bibitem{Stconj}J. Steenbrink,
\textit{The spectrum of hypersurface singularities}, Actes
du Colloque de Th\'eorie de Hodge (Luminy, 1987).
Ast\'erisque No. 179-180, (1989) \textbf{11}, 163--184.

\bibitem{sumi}H. Sumihiro,
\textit{Equivariant completion II},  J. Math. Kyoto Univ.,
(1975) \textbf{15}, 573--605.



\bibitem{Va}A. Varchenko, \textit{Asymptotic Hodge structure in the vanishing
cohomology}, Math. USSR Izvestija \textbf{18} (1982),
469--512.


\bibitem{v1}O. Villamayor, \textit{Constructiveness of Hironaka's resolution},
Ann. Sci. {\'E}c. Norm. Sup. Paris, \textbf{22} (1989), 1--32.


\bibitem{v2}O. Villamayor, \textit{Patching local uniformizations},
Ann. Sci. {\'E}c. Norm. Sup. Paris, \textbf{25} (1992), 629--677.

\end{thebibliography}

\end{document}